\documentclass[12pt, twoside]{article}
\usepackage{amsmath,amsthm,amssymb}
\usepackage{times}
\usepackage{enumerate}
\usepackage{latexsym,mathrsfs,bbm,amsbsy,mathdots}
\usepackage[mathx]{mathabx}
\allowdisplaybreaks
\def\titlerunning#1{\gdef\titrun{#1}}
\makeatletter
\def\author#1{\gdef\autrun{\def\and{\unskip, }#1}\gdef\@author{#1}}
\def\address#1{{\def\and{\\\hspace*{18pt}}\renewcommand{\thefootnote}{}%
\footnote {#1}}%
\markboth{\autrun}{\titrun}}
\makeatother
\def\email#1{e-mail: #1}
\def\subjclass#1{{\renewcommand{\thefootnote}{}%
\footnote{\emph{Mathematics Subject Classification (2010):} #1}}}

\newtheorem{theorem}{Theorem}[section]
\newtheorem{corollary}[theorem]{Corollary}
\newtheorem{lemma}[theorem]{Lemma}
\theoremstyle{definition}
\newtheorem{definition}[theorem]{Definition}
\newtheorem{remark}[theorem]{Remark}

\numberwithin{equation}{section}
\def\noteq{\operatorname{=\hspace{-3.8mm}/\ }}
\def\nequiv{\operatorname{\equiv\hspace{-3.7mm}/\ }}

\begin{document}


\baselineskip=17pt


\titlerunning{Weyl-Schr{\"o}dinger representations}

\title{Weyl-Schr\"odinger representations of Heisenberg groups in infinite dimensions}

\author{Oleh Lopushansky}

\date{}

\maketitle
\vspace{-2cm}

\address{Institute of Mathematics,
University of Rzesz\'ow,  \email{ovlopusz@ur.edu.pl}}



\begin{abstract}
We investigate the  group $\mathcal{H}_\mathbb{C}$ of complexified Heisenberg matrices with entries from an infinite-dimensional complex Hilbert space $H$.
Irreducible representations of the Weyl--Schr{\"o}dinger type on the space $L^2_\chi$ of quadratically integrable $\mathbb{C}$-valued functions are described.
Integrability is understood with respect to  the projective limit $\chi=\varprojlim\chi_i$ of probability Haar measures $\chi_i$ defined on groups of unitary $i\times i$-matrices $U(i)$. The measure $\chi$ is invariant under the infinite-dimensional group $U(\infty)=\bigcup U(i)$ and satisfies the abstract Kolmogorov consistency conditions. The space $L^2_\chi$ is generated by Schur polynomials on Paley--Wiener maps. The Fourier-image of $L^2_\chi$ coincides with the Hardy space ${H}^2_\beta$ of Hilbert--Schmidt analytic functions on $H$ generated by the correspondingly weighted Fock space $\Gamma_\beta(H)$.
An application to heat equation over $\mathcal{H}_\mathbb{C}$ is considered.
\medskip


\end{abstract}


\section{Introduction}\label{1}

An aim of this work is to investigate irreducible Weyl--Schr{\"o}dinger representations of the complexified Heisenberg group $\mathcal{H}_\mathbbm{C}$ (see \cite[n.9]{Lopushansky2017}),  consisting  of matrix elements  $X(a,b,t)$ with any $a,b\in{H}$ and $t\in\mathbb{C}$ such that
\begin{equation}\label{HC}
X(a,b,t)=\begin{bmatrix}
                      1 & a&t \\
                      0 &\mathbbm{1}&b \\
                      0&0&1
                    \end{bmatrix},\quad
X(a,b,t)\cdot X(a',b',t')=\begin{bmatrix}
                      1&a+a'&
                      t+t'+\langle{a}\mid{b'}\rangle\\
                      0 &\mathbbm{1}&b+ b'\\
                      0&0&1
                    \end{bmatrix}
\end{equation}
where $H$ is an infinite-dimensional  complex Hilbert space and $\mathbbm{1}$ is its identity map.

The group $\mathcal{H}_\mathbbm{C}$ has the unit $X(0,0,0)$ and inverse  elements of the form $X(a,b,t)^{-1}=X\left(-a,-b,-t+\langle{a}\mid{b}\rangle\right)$.

In what follows, we consider the infinite-dimensional unitary group $U(\infty)=\bigcup{U}(i)$, containing all subgroups
$U(i)$ of unitary $i\times i$-matrices, which acts irreducibly on a complex Hilbert space $\left\{H,\langle\cdot\mid\cdot\rangle\right\}$ with an orthonornal basis $(\mathfrak{e}_i)_{i\in\mathbb{N}}$.

To find the desired representation, we use the space $L^2_\chi$ of $\mathbb{C}$-valued functions that are quadratically integrable with respect
to the probability measure $\chi$. Wherein, according to our assumption $\chi$ has a structure of the projective limit $\chi=\varprojlim\chi_i$ of probability Haar's measures $\chi_i$ on $U(i)$, satisfying the Kolmogorov consistency conditions in an abstract Bochner's  formulation  (see \cite{Okada1977,Rao1971}).

In  \cite{zbMATH01820924,zbMATH02037136} it was shown that the projective limit $\chi=\varprojlim\chi_i$ is well defined over the projective limit
$\mathfrak{U}=\varprojlim U(i)$ with respect to the Liv\v{s}ic transforms $\pi_i^{i+1}\colon {U(i+1)\to U(i)}$ such that $\chi_i=\pi_i^{i+1}(\chi_{i+1})$.
In this paper, we prove that  for such $\chi$ each function from $L^2_\chi$ admit a superposition (linearization in the sense of \cite{CarZald04}) on Paley–Wiener maps
associated with $U(\infty)$. As a result,  it is shown that Schur polynomials  form an orthonormal basis in $L^2_\chi$ and the
Fourier-image of $L^2_\chi$ consists of Hilbert-Schmidt analytic functions on $H$.

 Note also that projective limits of probability measures over various  infinite-dimensional manifolds with similar properties
were investigated in \cite{Yamasaki73,Yamasaki85,Pickrell}.

If instead of the unitary group $U (\infty)$ we take the infinite-dimensional linear space with a Gaussian measure $\gamma$, a similar
construction of the appropriate space $L^2_\gamma$ can be found in the well-known works \cite{Bargmann1961}.
In this case,  the Fourier-image of $L^2_\gamma$ coincides with the Segal--Bargmann space of entire analytic functions
over which the Schr\"{o}dinger type irreducible representations of Heisenberg groups are well defined.
In the present paper, we change $\gamma$ by the unitarily-invariant projective limit  $\chi=\varprojlim\chi_i$
and, as a result, we obtain another irreducible representation, called to be the Weyl-Schr\"{o}dinger type.

Infinite-dimensional Heisenberg groups over $\mathbb{R}$ was considered in \cite{Neeb2000} by using the reproducing kernel Hilbert spaces.
The Schr\"{o}dinger representation of such groups using Gaussian measures over a real Hilbert space  was described in \cite{Beltita16}.
Since the group $\mathcal{H}_\mathbbm{C}$ in the case of matrix  entries $a,b,t\in\mathbb{R}$  coincides with the classical Heisenberg group
over $\mathbb{R}$ (see, e.g. \cite{Hall2013}), the results of the present paper can be considered as a complexification of previous studies.
The Weyl-Schr\"{o}dinger representation obtained here is not equivalent to that was described earlier.

Further, let us briefly describe the main results.
Consider the following mapping $\phi\colon H\ni{h}\longmapsto \phi_h\in L^2_\chi$ defined by Paley--Wiener maps
\begin{equation}\label{trace}
{\phi}_h(\mathfrak{u}):=
\sum{\phi}_i(\mathfrak{u})\,\mathfrak{e}_i^*(h) \quad \text{with}\quad
{\phi}_i(\mathfrak{u}):=\left\langle u_i(\mathfrak{e}_i)\mid\mathfrak{e}_i\right\rangle,
\quad u_i=\pi_i(\mathfrak{u}),
\end{equation}
where $\mathfrak{e}^*_i(\cdot):=\left\langle\cdot\mid\mathfrak{e}_i\right\rangle$ and  the
projections $\pi_i\colon\mathfrak{U}\ni\mathfrak{u} \to u_i\in U(i)$ are uniquely defined by $\pi_i^{i+1}$.
Every function $\phi_h$ of variable ${\mathfrak{u}\in\mathfrak{U}}$
satisfies the equality (Cor.\ref{ggauss})
\begin{equation*}\label{gau}
\int\exp\big\{\mathop{\mathrm{Re}}\phi_h\big\}\,d\chi=
 \exp\Big\{\frac{1}{4}\|h\|^2\Big\},\quad{h\in H}.
\end{equation*}

The space $L^2_\chi $ can be generated by two orthonormal bases, consisting of
Schur polynomials and power polynomials  of variables  ${\phi}_\imath=
\left({\phi}_{\imath_1},\ldots,{\phi}_{\imath_\eta}\right)$, respectively,
\begin{equation}\label{schur0}
s^\lambda_\imath(\mathfrak{u}):=\frac{\mathop{\mathrm{det}}
\big[\phi_{\imath_i}^{\lambda_j+\eta-j}(\mathfrak{u})\big]_{1\le i,j\le\eta}}
{\prod_{1\le i<j\le\eta}[\phi_{\imath_i}(\mathfrak{u})-\phi_{\imath_j}(\mathfrak{u})]}
\quad\text{and}\quad
{\phi}^\lambda_\imath:=
{\phi}^{\lambda_1}_{\imath_1}\ldots{\phi}^{\lambda_\eta}_{\imath_\eta}.
\end{equation}
These bases are indexed by tabloids $\imath^\lambda$ with strictly ordered $\imath=({\imath_1},\ldots,{\imath_\eta})\in\mathbb{N}^\eta$
where $\lambda=(\lambda_1,\ldots,\lambda_\eta)\in\mathbb{N}^\eta$ is a partition of $n\in\mathbb{N}$ and $\eta=\eta(\lambda)$ stands for the length of $\lambda$.
Then we write briefly ${\imath^\lambda\vdash n}$. The orthogonal expansion $L_\chi^2 =\bigoplus{L}_\chi^{2,n}$ holds (Thm~\ref{Schur1})
where ${L}_\chi^{2,n}$ are formed by $n$-homogeneous polynomials ${\phi}^\lambda_\imath$, normed as follows
\[
\|{\phi}^\lambda_\imath\|_\chi^2=\int |{\phi}^\lambda_\imath|^2d\chi=\beta_\lambda
\lambda!,\qquad\beta_\lambda:=\dfrac{(\eta-1)!}{(\eta-1+n)!},  \quad\lambda!:=\lambda_1!\ldots\lambda_\eta!.
\]

It is also shown  that the surjective linear isometry $\varPsi\colon{H}^2_\beta\ni\psi^*_f\longmapsto{f}\in{L}^2_\chi$
holds (Lem.~\ref{pbasis}), where ${H}^2_\beta=\sum{P}_\beta^n(H)$ means the Hardy space of  entire analytic functions  $\psi^*_f(h)$ of variable $h\in H$ and
${P}_\beta^n(H)$ is generated by the $n$-homogeneous Hilbert--Schmidt polynomials $\mathfrak{e}^{*\lambda}_\imath:=
\mathfrak{e}^{*\lambda_1}_{\imath_1}\ldots\mathfrak{e}^{*\lambda_\eta}_{\imath_\eta}$,  normed as
$\|\mathfrak{e}^{*\lambda}_\imath\|_{{H}^2_\beta}= \big(\beta_\lambda{\lambda!}\big)^{1/2}$.

If the basis of symmetric tensor elements $\mathfrak{e}^{\odot\lambda}_\imath:=\mathfrak{e}^{\otimes\lambda_1}_{\imath_1}\odot\ldots\odot
\mathfrak{e}^{\otimes\lambda_\eta}_{\imath_\eta}$ (associated with $\mathfrak{e}^{*\lambda}_\imath$)
in the correspondingly weighted Fock space $\Gamma_\beta(H)$
is normed as $\|\mathfrak{e}^{\odot\lambda}_\imath\|_{\Gamma_\beta}=
\|\mathfrak{e}^{*\lambda}_\imath\|_{{H}^2_\beta}$ then each function $f\in{L}^2_\chi$ admits the superposition
\[
f=\varPsi\circ\psi^*_f,\qquad\psi^*_f(h)
=\sum_{n\ge0}\frac{1}{n!}\sum_{\imath^\lambda\vdash n}\frac{n!}{\lambda!}
\mathfrak{e}^{*\lambda}_\imath(h)
\big\langle\mathfrak{e}^{\odot\lambda}_\imath\mid\psi_f\big\rangle_{\Gamma_\beta},
\quad{h\in{H}},
\]
where the Taylor expansion on the right-hand side
of any analytic function $\psi^*_f\in {H}^2_\beta$ on $H$ is uniquely determined by the corresponding element
$\psi_f\in\Gamma_\beta(H)$.

Our further goal is to analyze  the inverse isomorphism $\varPsi^{-1}$ which can be described by the Fourier transform under the measure $\chi$
in following way
\[
\hat{f}(h)=\int \exp(\bar\phi_h)f\,d\chi\quad\text{where}\quad
F=\varPsi^{-1}\colon{L}^2_\chi \ni{f}\longmapsto\hat{f}:=\psi^*_f\in{H}^2_\beta.
\]
The Fourier transform $F$ acts isometrically on the Hardy space of analytic functions ${H}^2_\beta$ (Thm~\ref{laplace1}).
So, $F$ acts as an analytic extension of the mapping $\phi$.

Applying the superposition with $\varPsi$, we describe two different representations
of the additive group  $(H,+)$ over  $L^2_\chi $  defined by shift and multiplicative groups (Lem.~\ref{comm}). Using this
we show (in Thm~\ref{heis}) that an irreducible representation of
the Heisenberg group $\mathcal{H}_\mathbb{C}$
can be realized on  ${L}^2_\chi $ in the Weyl--Schr\"{o}dinger form
\[
X(a,b,z)\longmapsto\exp(z){W}^\dagger(a,b),\quad{W}^\dagger(a,b):=
\exp\Big\{\frac{1}{2}\langle{a}\mid{b}\rangle\Big\}T^\dagger_bM^\dagger_{a^*}
\]
for all $a,b\in H$ and $z\in\mathbb{C}$, where $T_b^\dagger$ and $M_{a^*}^\dagger$ are defined by shift and multiplicative groups, respectively.
It is also proved that the Weyl system ${W}^\dagger(a,b)$
has  the densely-defined generator $\mathfrak{p}^\dagger_{a,b}:=\partial_{b}^\dagger+\bar\phi_{a}$ which satisfies the
commutation relation
\[
{W}^\dagger(a,b){W}^\dagger(a',b')
=\exp\left\{ - \big[\mathfrak{p}_{a,b}^\dagger,\mathfrak{p}_{a',b'}^\dagger\big]\right\}
{W}^\dagger(a',b'){W}^\dagger(a,b)
\]
where the groups $M_{a^*}^\dagger$ and $T_b^\dagger$ are generated  by $\bar\phi_a$ and $\partial_b^\dagger$, respectively.

Applying the Weyl--Schr\"{o}dinger representation to the associated with $\mathcal{H}_\mathbb{C}$ heat equation, we prove (Thm~\ref{schrod})
that the following Cauchy problem with $\partial_i^\dagger:=\partial_{\mathfrak{e}_i}^\dagger$,
\[
\frac{dw(r)}{dr}=-\sum\big(\partial_i^\dagger+\bar\phi_i\big)^2w(r),
\quad w(0)=f,\quad r>0,
\]
has the unique solution $w(r)=\mathfrak{G}^\dagger_rf$ for any function $f$ from a finite sum  $\bigoplus{L}^{2,n}_\chi$,
where the $1$-parameter Gaussian semigroup $\mathfrak{G}_r^\dagger$ has the form
\begin{align*}
\mathfrak{G}^\dagger_rf&=\frac{1}{\sqrt{4\pi r}}\int_{c_0}
\exp\Big\{-\frac{\|\tau\|_{w_0}^2}{4r}\Big\}{W}^\dagger_\tau{f}\,d\mathfrak{w}(\tau),\\
{W}^\dagger_\tau{f}&:=\lim_{n\to\infty}
\exp\Big\{-\frac{\|{p_n^\sim(\tau)}\|_{w_0}^2}{2}\Big\}
\prod_{i=1}^{n} T^\dagger_{\mathbbm{i}\tau_i\mathfrak{e}_i}
M^\dagger_{-\mathbbm{i}\tau_i\mathfrak{e}_i^*}.
\end{align*}
Here $\tau=(\tau_i)$ belongs to the abstract Wiener space
$\{w_0,\|\cdot\|_{w_0}\}$ defined by the injections $l_2\looparrowright{w}_0\looparrowright{c}_0$ of real Banach spaces
and endowed with the Wiener measure $\mathfrak{w}$  in according to the known Gross' theorem \cite{Gross65},
whereas the sequence of projectors $(p^\sim_n)$ onto $\mathbb{R}^n$ is convergent  to the identity map on $w_0$.

Finally, note that this work is a continuation of previous publications \cite{Lopushansky2016,Lopushansky2017}.
The novelty results from the observation that the system of Schur polynomials with variables on Paley--Wiener maps
form an orthonormal basis in $L^2_\chi $. This allowed us to investigate irreducible
Weyl--Schr\"{o}dinger representations and  Weyl systems of the Heisenberg group $\mathcal{H}_\mathbb{C}$ on the whole space  $L^2_\chi$.

\section{Invariant probability measure}\label{2}

Consider the  unitary group $U(\infty)=\bigcup U(m)$ with $m\in\mathbb{N}_0=\mathbb{N}\cup\{0\}$,
$ \mathbbm{1}=U(0)$, irreducibly acting on a separable Hilbert space ${H}$, where subgroups $U(m)$ are identified  with ranges of injections
$U(m)\ni u_m\longmapsto  {\begin{bmatrix}
                      u_m & 0\\
                      0 &\mathbbm{1}\\
                    \end{bmatrix}\in U(\infty)}$. Following to \cite{zbMATH01820924}, \cite{zbMATH02037136},
we use the Liv\v{s}ic transforms $\pi^{m+1}_m\colon U(m+1)\to U(m)$ of the form
\begin{align}\label{projective}
&\pi^{m+1}_m\colon{u_{m+1}}:=\begin{bmatrix}
                      z_m &a \\
                      b & t \\
                    \end{bmatrix}\longmapsto u_m:=\left\{\begin{array}{clc}
                      z_m-[a(1+t)^{-1}b]&:& t\noteq -1 \\
                      z_m&:& t=-1 \\
                    \end{array}\right.
\end{align}
with $z_m\in U(m)$ defined by excluding $x_1=y_1\in\mathbb{C}$  from $\begin{bmatrix}
                      y_m \\
                      y_1\\
                    \end{bmatrix}=\begin{bmatrix}
                      z_m & a\\
                      -b &-t\\
                    \end{bmatrix}\begin{bmatrix}
                      x_m \\
                      x_1\\
                    \end{bmatrix}$
for $x_m,y_m\in\mathbb{C}^m$ and $a,b\in\mathbb{C}$ \cite[Lem.~3.1]{zbMATH02037136}.
It is surjective (not continuous) Borel mapping  \cite[Lem.~3.11]{zbMATH02037136}.

The projective limit $\mathfrak{U}:=\varprojlim U(m)$ under $\pi^{m+1}_m$
has  surjective Borel (not group homomorphisms) projections
\[\pi_m\colon{\mathfrak{U}\ni\mathfrak{u}\longmapsto u_m\in U(m)}\quad
\text{such that}\quad\pi_m={\pi^{m+1}_m\circ\pi_{m+1}}.\]
Their elements ${\mathfrak{u}\in\mathfrak{U}}$ are called the {\em virtual unitary matrices}.
The right action
 \[
 \mathfrak{U}\ni\mathfrak{u}\longmapsto
\mathfrak{u}.g\in\mathfrak{U}\quad\text{with}\quad g={(v,w)\in{U(\infty)\times U(\infty)}}
\]
is defined to be $\pi_m(\mathfrak{u}.g)=w^{-1}\pi_m(\mathfrak{u})v$, where $m$ is large enough that $v,w\in{U}(m)$.
On $\mathfrak{U}$ the involution $\mathfrak{u}\mapsto\mathfrak{u}^\star=(u_k^\star)$ is well defined, where
$u_k^\star=u_k^{-1}$ is adjoint to ${u_k\in U(k)}$. Thus, $[\pi_m(\mathfrak{u}.g)]^\star=\pi_m(\mathfrak{u}^\star.g^\star)$
for all $g^\star={(w^\star,v^\star)}\in{U(\infty)\times U(\infty)}$.

There exists the dense embedding $U(\infty)\looparrowright\mathfrak{U}$ (see  \cite[n.4]{zbMATH02037136}) which
assigns the stabilized sequence $\mathfrak{u}=(u_k)$ to each ${u_m\in U(m)}$  such that
\begin{equation}\label{stab}
\begin{split}
U(m)&\ni u_m\longmapsto(u_k)\in\mathfrak{U},\\
u_k&=\left\{\begin{array}{cl}
\pi^m_k(u_m)=(\pi^{k+1}_k\circ\ldots\circ \pi^m_{m-1})(u_m)&:k<m,\\
                      u_m &: k\ge m.
\end{array}\right.\end{split}
\end{equation}

We always assume that the group $U(m)$ is endowed with the probability Haar measure $\chi_m$.
Using the Kolmogorov consistency theorem (see, e.g.
\cite[Lem.4.8]{zbMATH02037136},  \cite[Thm 2.2]{Rao1971}, \cite[Cor.4.2]{Tomas2006}),
we determine the probability measure
on $\mathfrak{U}$ to be the projective limit \[\chi:=\varprojlim\chi_m\quad\text{under}\quad
\chi_m=\pi_m^{m+1}(\chi_{m+1})
\]
where $\pi_m^{m+1}(\chi_{m+1})$ means an image-measure and ${\chi_0=1}$.
 As is known   \cite[Thm 2.5]{Tomas2006}, the measure $\chi$ is  Radon.
 We now describe the necessary properties of $\chi$.

Consider the Hilbert space $L^2_\chi$ of functions $f\colon\mathfrak{U}\to\mathbb{C}$ with the following norm and inner product
\[\|f\|_\chi=\langle f\mid f\rangle_\chi^{1/2},\quad
\langle f_1\mid f_2\rangle_\chi:=\int f_1\bar f_2\,d\chi.\]
Let $L_\chi^\infty$ be the space of  $\chi$-essentially bounded functions $f\colon\mathfrak{U}\to\mathbb{C}$  with the norm
$\|f\|_\infty=\mathop{\rm ess\,sup}_{\mathfrak{u}\in\mathfrak{U}}|f(\mathfrak{u})|$.
The embedding $L^\infty_\chi\looparrowright L^2_\chi$  holds and $\|f\|_\chi\le\|f\|_\infty$.

\begin{lemma}\label{proh0}
For any $f\in{L}_\chi^\infty$ there exists the limit
\begin{equation}\label{occurs0}
\int f\,d\chi=\lim\int f\,d(\chi_m\circ\pi_m)=\lim\int (f\circ\pi_m^{-1})\,d\chi_m.
\end{equation}
Moreover, the measure $\chi$ is invariant under the right action, which means that
\begin{align}\label{Inv}
\int f(\mathfrak{u}.g)\,d\chi(\mathfrak{u})&=\int f(\mathfrak{u})\,d\chi(\mathfrak{u}),
\quad g\in U(\infty)\times U(\infty),\\\label{inv0}
\int f\,d\chi&=\int d\chi(\mathfrak{u})
\int f(\mathfrak{u}.g)\,d(\chi_m\otimes\chi_m)(g).
\end{align}
\end{lemma}

\begin{proof}
The sequence $\{(\chi_m\circ\pi_m)(\mathcal{K})\}$ is decreasing
for any compact set $\mathcal{K}$ in $\mathfrak{U}$, since $\pi_m={\pi^{m+1}_m\circ\pi_{m+1}}$ yields
$\pi_{m+1}(\mathcal{K})\subseteq(\pi^{m+1}_m)^{-1}\left[\pi_m(\mathcal{K})\right]$. It follows
\begin{equation}\label{decreas}\begin{split}
(\chi_m\circ\pi_m)(\mathcal{K})&=\pi^{m+1}_m(\chi_{m+1})\left[\pi_m(\mathcal{K})\right]\\
&=\chi_{m+1}\left[(\pi^{m+1}_m)^{-1}[\pi_m(\mathcal{K})]\right]
\ge(\chi_{m+1}\circ\pi_{m+1})(\mathcal{K}).
\end{split}\end{equation}
This ensures that the necessary and sufficient conditions of the Prokhorov theorem \cite[Thm IX.52]{BourbakiINTII}
and its modification from \cite[Thm 4.2]{Tomas2006} are satisfied.

Indeed, let $\check{U}(m)\subset U(m)$ be the set of matrices with no eigenvalue $\{-1\}$ for ${m\ge1}$. As is known \cite[n.3]{zbMATH02037136}, $\check{U}(m)$ is open in $U(m)$ and  ${\chi_m(U(m)\setminus\check{U}(m))=0}$. In virtue of \cite[Lem. 3.11]{zbMATH02037136}
the restrictions $\pi^{m+1}_m\colon\check{U}(m+1)\rightarrow\check{U}(m)$ are continuous and surjective. The projective limit $\varprojlim\check{U}(m)$ under these restrictions has continuous surjective projections $\pi_m\colon\varprojlim\check{U}(m)\rightarrow\check{U}(m)$. Restrict $\chi_m$ to $\check{U}(m)$.
By \cite[Thm 6]{Tomas2006}, a probability measure $\check{\chi}$ satisfying conditions $\pi_m(\check{\chi})=\chi_m|_{\check{U}(m)}$ is well defined
iff for every $\varepsilon>0$ there exists a compact set $\mathcal{K}\subset\varprojlim\check{U}(m)$ such that
\[
{(\chi_m\circ\pi_m)(\mathcal{K})}\ge{1-\varepsilon}\quad\text{for all}\quad
{m\in\mathbb{N}}.
\]
Then by the Prokhorov theorem $\check{\chi}$ is uniquely determined as
\begin{equation}\label{Prokh}
\check{\chi}(\mathcal{K})=\inf(\chi_m\circ\pi_m)(\mathcal{K})\quad\text{for all}\quad\mathcal{K}\subset\varprojlim\check{U}(m).
\end{equation}

Let $\varepsilon>0$ and $K_1\subset\check{U}(1)$ be a compact set such that $\chi_1(K_1)>1-\varepsilon$. Let
a compact sets $K_m\subset\check{U}(m)$ be defined inductively   such that
\[
\pi^{m+1}_m(K_{m+1})\subset K_m\quad\text{and}\quad  \chi_{m+1}(K_{m+1})>1-\varepsilon\quad\text{for all}\quad
 {m\ge1}.
 \]
 Assume that $K_1,\ldots,K_m$ are constructed.  Since $\chi_m=\pi^{m+1}_m(\chi_{m+1})$, we get
\[
\chi_m(K_m)=\chi_{m+1}[(\pi^{m+1}_m)^{-1}(K_m)]>1-\varepsilon.
\]
By regularity of $\chi_{m+1}|_{\check{U}(m)}$, there exists a compact set
\[
K_{m+1}\subset(\pi^{m+1}_m)^{-1}(K_m)\quad\text{such that}\quad\chi_{m+1}(K_{m+1})>1-\varepsilon.
\]
The induction is complete. Then $\mathcal{K}=\varprojlim K_m$
with $K_0=\mathbbm{1}$ is compact. By virtue of \eqref{Prokh}, we have $\check\chi(\mathcal{K})\ge1-\varepsilon$. Hence, the projective limit
$\check\chi=\varprojlim\chi_m|_{\check{U}(m)}$ is well defined on $\varprojlim\check{U}(j)$ by the Prokhorov criterion.

The measure $\check{\chi}$ can be extended to $\varprojlim{U}(m)\setminus\varprojlim\check{U}(m)$ as zero,
since each $\chi_m$ is zero on ${U(m)\setminus\check{U}(m)}$. The uniqueness of the projective limits yields $\check{\chi}=\chi$.
So,  $\chi=\varprojlim\chi_m$ is also well defined and by  \eqref{decreas} and \eqref{Prokh} we get
\[
\chi(\mathcal{K})=\inf(\chi_m\circ\pi_m)(\mathcal{K})=\lim(\chi_m\circ\pi_m)(\mathcal{K})\quad\text{for all compact}\quad
\mathcal{K}\subset\mathfrak{U}.
\]
By the known Portmanteau theorem \cite[Thm 13.16]{Klenke2008} it  follows that the limit \eqref{occurs0} exists. Whereas,
the property \eqref{Inv} is a consequence of the equalities
\[
\chi(\mathcal{K}.g)=\lim\chi_m(K_m.g)=\lim\chi_m(K_m)=\chi(\mathcal{K})
\]
for all $g=(v,w)\in{U(\infty)\times U(\infty)}$ where $m$ is large enough that $v,w\in U(m)$.

Finally, the function $(\mathfrak{u},g)\mapsto f(\mathfrak{u}.g)$ with any
${f\in{L}_\chi^\infty}$  is integrable over ${\mathfrak{U}\times U(m)\times U(m)}$, hence
\[
\int\,d\chi(\mathfrak{u})\int f(\mathfrak{u}.g)\,d(\chi_m\otimes\chi_m)(g)
=\int\,d(\chi_m\otimes\chi_m)(g)\int f(\mathfrak{u}.g)\,d\chi(\mathfrak{u})
\]
by the Fubini theorem. It yields  \eqref{inv0}
since the internal integral on the right-hand side is independent  of $g$ by  \eqref{Inv} and
${\int\,d(\chi_m\otimes\chi_m)(g)=1}.$
The proof  is complete.
\end{proof}

We now note  the concentration property of Haar measures sequence $(\chi_m)$ satisfying the Kolmogorov conditions
$\chi_m=\pi_m^{m+1}(\chi_{m+1})$ if each group $U(m)$ is endowed with the normalized Hilbert--Schmidt metric
\[d_{HS}(u,v)=\sqrt{m^{-1}\mathop{\sf tr}|u-v|_{HS}}\quad \text{where}\quad |u-v|_{HS}=\sqrt{(u-v)^\star(u-v)}.\]

As is well known (see \cite{Milman83,Voiculescu1991}),  $(U(m),d_{HB},\chi_m)$ is a L\'{e}vy family. Namely,
the following  sequence of isoperimetric constants dependent on $\varepsilon>0$
\[\alpha(U(m),\varepsilon)=1-\inf\big\{\chi_m[(\Omega_m)_\varepsilon]\colon \Omega_m \text{ be Borel set in } U(m), \chi_m(\Omega_m)>1/2\big\}\]
with $(\Omega_m)_\varepsilon=\left\{u_m\in U(m)\colon d_{HS}\left(u_m,\Omega_m\right)<\varepsilon\right\}$
is such that \[\alpha(U(m),\varepsilon)\to0\quad \text{as}\quad{m\to\infty}.\]
Taking into account the Lemma~\ref{proh0}, we can formulate the following conclusion.

\begin{corollary}\label{Levi}
For any Borel set $\Omega_\varepsilon=\varprojlim (\Omega_m)_\varepsilon$ with $\chi_m(\Omega_m)>1/2$
in the projective limit $\mathfrak{U}=\varprojlim U(m)$  the equality
\[{\chi(\Omega_\varepsilon)=\lim_{m\to\infty}\chi_m\left[(\Omega_m)_\varepsilon\right]=1}
\]
 holds. Consequently, all Borel sets $\mathfrak{U}\setminus\Omega_\varepsilon$ with ${\chi_m(\Omega_m)>1/2}$ and any ${\varepsilon>0}$
are $\chi$-measure zero, i.e., the measure $\chi=\varprojlim\chi_m$ is concentrated outside these sets.
\end{corollary}

\section{Polynomials on Paley--Wiener maps}\label{3}

Let $\mathscr{I}_\eta:=\big\{\imath=\big({\imath_1},\ldots,{\imath_\eta}\big)\in\mathbb{N}^\eta
\colon\imath_1<\imath_2<\ldots<\imath_\eta\big\}$ be an integer alphabet of
length $\eta$ and $\mathscr{I}=\bigcup\mathscr{I}_\eta$.
Let $\lambda=(\lambda_1,\ldots,\lambda_\eta)\in\mathbb{N}^\eta$ with
${\lambda_1\ge\lambda_2\ge\ldots\ge\lambda_\eta}$ be a partition of  an $n$-letter  word
$\imath^\lambda=\big\{\Box_{ij}\colon {1\le i\le\eta}, j=1,\ldots,\lambda_i\big\}$
with  ${\imath\in\mathscr{I}_\eta}$.
A Young $\lambda$-tableau with a partition $\lambda$
is a result of filling the word  $\imath^\lambda$ onto the  matrix
$[\imath^\lambda]=\begin{array}{cccc}
    \Box_{11} & \dots & \dots &\Box_{1\lambda_1}\\
    \vdots & \vdots & \iddots &    \\
    \Box_{\eta1}  & \dots &\Box_{\eta\lambda_\eta} &
\end{array}\!\!\!\!\!$
with  $n$ nonzero entries in some way without repetitions. So,
each $\lambda$-tableau $[\imath^\lambda]$ can be identified with a bijection
$[\imath^\lambda]\to\imath^\lambda$.
The conjugate partition $\lambda^\intercal$ corresponds to the transpose matrix $[\imath^\lambda]^\intercal$.

A Young tableau $[\imath^\lambda]$ is called {\em standard} ({\em semistandard}\,) if its entries are strictly (weakly)
ordered along each row and strictly ordered down each column. Let $\mathbb{Y}$ denote all Young tabloids $[\imath^\lambda]$
 and $\mathbb{Y}_n$ be its subset such that $\imath^\lambda\vdash n$.
Assume that $\mathbb{Y}_0=\left\{\emptyset\in\mathbb{Y}\colon|\emptyset|=0\right\}$
and $\eta(\emptyset)=0$.

As before, $\left\{H,\langle\cdot\mid\cdot\rangle\right\}$ is a separable complex Hilbert space with an orthonormal basis
  $\left\{\mathfrak{e}_i\colon i\in\mathbb{N}\right\}$ and ${\|\cdot\|={\langle\cdot\mid\cdot\rangle^{1/2}}}$.
For  its adjoint space $H^*$ the conjugate-linear isometry ${*\colon H^*\rightarrow H^{**}=H}$ is defined via
  $a^*(h)={\langle h\mid a\rangle}$ for all ${a,h\in H}$. The Fourier expansion
$h=\sum\mathfrak{e}^*_i(h)\mathfrak{e}_i$ with
$\mathfrak{e}^*_i(h):={\langle h\mid\mathfrak{e}_i\rangle}$  holds.
The tensor power ${H}^{\otimes n}\!,$ spanned  by elements $\psi_n={h_1\otimes\ldots\otimes h_n}$  with ${h_i\in{H}}$  ${(i=1,\ldots,n)}$,
is endowed with the norm
$\|\psi_n\|={\left\langle\psi_n\mid\psi_n\right\rangle}^{1/2}$
where  ${\left\langle\psi_n\mid\psi_n'\right\rangle}:=
{\langle h_1\mid h'_1\rangle\ldots\langle h_n\mid h'_n\rangle}$.

Let $S_n$ be the group of $n$-elements permutations $\sigma(\psi_n):={h_{\sigma(1)}\otimes\ldots\otimes{h}_{\sigma(n)}}$.
An orthogonal basis in ${H}^{\otimes n}$ is formed by elements
$\sigma(\mathfrak{e}^{\otimes\lambda_1}_{\imath_1}\otimes\ldots\otimes
\mathfrak{e}^{\otimes\lambda_\eta}_{\imath_\eta})$ with ${\imath^\lambda\vdash n}$ and $\eta=\eta(\lambda)$,
additionally indexed by all $\sigma\in S_n$.
The symmetric tensor power $H^{\odot n}\subset H^{\otimes n}$ is defined to be a range of the orthogonal projector
${\mathcal{S}_n\colon{H}^{\otimes n}\ni\psi_n\longmapsto{h_1\odot\ldots\odot h_n}:=(n!)^{-1}{\sum}_{\sigma\in S_n}\sigma(\psi_n)}$.
We assume that  ${H}^{\otimes n}$ is completed and that ${H}^{\otimes 0}=\mathbb{C}$.
Let $\psi_n:=h^{\otimes n}$ for ${h=h_i}$. The embedding
$\left\{h^{\otimes n}\colon {h}\in{H}\right\}\subset{H}^{\odot n}$ is total  by  the polarization formula \cite[n.1.5]{Floret97}
\begin{equation}\label{polarization}
{h}_1\odot\ldots \odot{h}_n=\frac{1}{2^nn!}
\sum_{\theta_1,\ldots,\theta_n=\pm 1} \theta_1\dots \theta_n{h}^{\otimes n},\quad
{h}=\sum_{i=1}^n\theta_ih_i.
\end{equation}

Let $H_\eta\subset H$ be spanned by $\big\{\mathfrak{e}_{\imath_1},\ldots,\mathfrak{e}_{\imath_\eta}\big\}$.
We can uniquely assign to any semistandard tableau $[\imath^\lambda]$ with $\imath^\lambda\vdash n$
the element in $H_\eta^{\otimes n}$ for which there exists the permutation ${\sigma'\in{S}_n}$ such that
${\sigma'\big(\mathfrak{e}^{\otimes\lambda_1}_{\imath_1}\otimes\ldots\otimes
\mathfrak{e}^{\otimes\lambda_\eta}_{\imath_\eta}\big)}=
{\mathfrak{e}^{\otimes\lambda_1}_{\imath_1}\odot\ldots\odot\mathfrak{e}^{\otimes\lambda_\eta}_{\imath_\eta}\in H_\eta^{\odot n}}$.
Taking all $\imath\in\mathscr{I}$, we conclude that the system indexed by semistandard $\lambda$-tabloids
\[
\begin{split}
\mathfrak{e}^{\mathbb{Y}_n}&=\left\{\mathfrak{e}^{\odot\lambda}_\imath
:=\mathfrak{e}^{\otimes\lambda_1}_{\imath_1}\odot\ldots\odot
\mathfrak{e}^{\otimes\lambda_\eta}_{\imath_\eta}
\colon\imath^\lambda\vdash n, \ \lambda\in\mathbb{Y}_n, \ \imath\in\mathscr{I}\right\},
\quad{\mathfrak{e}^{\odot\emptyset}_\imath=1}\\
\text{where}&\quad
\langle\mathfrak{e}^{\odot\lambda}_\imath\mid
\mathfrak{e}^{\odot\lambda'}_{\imath'}\rangle
=\left\{
\begin{array}{clcl}
{\lambda!}/{n!}
&: \lambda=\lambda&\text{ and }&\imath=\imath' \\
  0 &: \lambda\noteq\lambda'&\text{ or }&\imath\noteq\imath'
\end{array}\right.
\end{split}
\]
forms an orthogonal basis in the symmetric tensor power $H_\eta^{\odot n}$.

The system $\big\{\mathfrak{e}_\imath^{\otimes\lambda}:=
\mathcal{S}_n\big(\mathfrak{e}^{\otimes\lambda_1}_{\imath_1}\otimes\ldots\otimes
\mathfrak{e}^{\otimes\lambda_\eta}_{\imath_\eta}\big)\colon\imath^\lambda\vdash n,
\ \lambda\in\mathbb{Y}_n, \ \imath\in\mathscr{I}\big\}$, additionally indexed by all $\sigma\in S_n$,
forms an orthonormal basis in the whole tensor power $H^{\otimes n}$.

As  usually, the {\em symmetric Fock space} is defined to be the Hilbertian orthogonal sum
$\Gamma(H)={\bigoplus}_{n\ge0} H^{\odot n}$  with the orthogonal basis
$\mathfrak{e}^{\mathbb{Y}}:={\bigcup\big\{\mathfrak{e}^{\mathbb{Y}_n}\colon{n}\in\mathbb{N}_0\big\}}$
of elements $\psi=\bigoplus\psi_n$ with $\psi_n\in H^{\odot n}$
endowed with the inner product and norm
\[
\langle\psi\mid\psi'\rangle_\Gamma=\sum{n!}\langle\psi_n\mid\psi_n'\rangle,\quad
\|\psi\|_\Gamma=\langle\psi\mid\psi\rangle^{1/2}_\Gamma.
\]

Note that by tensor multinomial theorem the Fourier expansion under $\mathfrak{e}^{\mathbb{Y}_n}$
\begin{equation}\label{Tayl0}
h^{\otimes n}=\sum_{\imath^\lambda\vdash n}
\frac{n!}{\lambda !}\mathfrak{e}^{\odot\lambda}_\imath\,\mathfrak{e}^{*\lambda}_\imath(h),\quad
\|h^{\otimes n}\|^2 =
\sum_{\imath^\lambda\vdash n}
\frac{n!}{\lambda !}|\mathfrak{e}^{*\lambda}_\imath(h)|^2,\quad \mathfrak{e}^{*\lambda}_\imath:=
\mathfrak{e}^{*\lambda_1}_{\imath_1}\ldots\mathfrak{e}^{*\lambda_\eta}_{\imath_\eta},
\end{equation}
holds in $H^{\odot n}$ for all $h\in H$. Consequently, the linearly independent, so-called,  coherent states ${\big\{\exp(h)\colon{h}\in{H}\big\}}$
in $\Gamma(H)$ have the expansion under the basis $\mathfrak{e}^\mathbb{Y}$
\begin{equation}\label{Tayl}
\exp(h):=\bigoplus_{n\ge0}\frac{h^{\otimes n}}{n!}
=\bigoplus_{n\ge0}\frac{1}{n!}\Big(\sum_{i\ge0}\mathfrak{e}_i\,\mathfrak{e}_i^*(h)\Big)^{\otimes n}=
\bigoplus_{n\ge0}\frac{1}{n!}\sum_{\imath^\lambda\vdash n}
\frac{n!}{\lambda!}\mathfrak{e}^{\odot\lambda}_\imath\,\mathfrak{e}^{*\lambda}_\imath(h)
\end{equation}
with $h^{\otimes 0}=1$, that is convergent, since
$\|\mathfrak{e}^{\odot\lambda}_\imath\|^2_\Gamma=
n!\|\mathfrak{e}^{\odot\lambda}_\imath\|^2 $ and
\begin{equation}\label{exp}
\begin{split}
\|\exp(h)\|^2_\Gamma&=\sum_{n\ge0}\frac{1}{n!}\sum_{\imath^\lambda\vdash n}
\Big(\frac{n!}{\lambda !}\Big)^2\|\mathfrak{e}^{\odot\lambda}_\imath\|^2
|\mathfrak{e}^{*\lambda}_\imath(h)|^2=\sum_{n\ge0}\frac{1}{n!}\sum_{\imath^\lambda\vdash n}
\frac{n!}{\lambda !}|\mathfrak{e}^{*\lambda}_\imath(h)|^2\\
&=\sum\frac{1}{n!}\left(\sum|\mathfrak{e}^*_i(h)|^2\right)^n=\sum\frac{1}{n!}\|h\|^{2n}=\exp\|h\|^2.
\end{split}
\end{equation}

\begin{definition}\label{PWm}
For any $h\in H$ and $\mathfrak{u}\in\mathfrak{U}$ the Paley--Wiener maps   are defined to be
\[
{\phi}_h(\mathfrak{u}):=
\sum{\phi}_i(\mathfrak{u})\,\mathfrak{e}_i^*(h) \quad \text{with}\quad
{\phi}_i(\mathfrak{u}):=\left\langle u_i(\mathfrak{e}_i)\mid\mathfrak{e}_i\right\rangle,
\quad u_i=\pi_i(\mathfrak{u})
\]
where  projections $\pi_i\colon\mathfrak{U}\ni\mathfrak{u} \to u_i\in U(i)$ are uniquely defined by $\pi_i^{i+1}$.
\end{definition}

These maps satisfy the orthogonal conditions $\phi_{\mathfrak{e}_i}=\phi_i$ and have the natural extension ${\phi}_{h^*}=\bar{{\phi}_h}$
onto the adjoint space $H^*$.

Note that, as in the case of linear spaces (see e.g. \cite[n.4.4]{Sengupta2014}, \cite{Stroock2010}), the Paley--Wiener  maps uniquely determine the embedding $\phi\colon H\ni{h}\longmapsto \phi_h\in L^2_\chi$.

For every ${h\in H}$ the ${l}_2$-valued function ${\phi}_h(\mathfrak{u})$ of variable ${\mathfrak{u}\in\mathfrak{U}}$
is well-defined,  since $(\mathfrak{e}_i^*(h))\in{l}_2$ and
$\left|\left\langle u_i(\mathfrak{e}_i)\mid\mathfrak{e}_i\right\rangle\right|\le1$.
We show that ${\phi}_h\in L^2_\chi$.   Assign for any partition $\lambda=(\lambda_1,\ldots,\lambda_\eta)\in\mathbb{N}^\eta$ of the weight $|\lambda|={\lambda_1+\ldots+\lambda_\eta}$
the constant
\begin{equation}\label{constant}
\beta_\lambda:=\dfrac{(\eta-1)!}{(\eta-1+|\lambda|)!}\le1,\quad\eta=\eta(\lambda).
\end{equation}

\begin{lemma}\label{proh2}
To every semistandard  tableau $[\imath^\lambda]$ one can uniquely assign the  function
\begin{equation}\label{FourF}
{\phi}^\lambda_\imath(\mathfrak{u}):=
{\phi}^{\lambda_1}_{\imath_1}(\mathfrak{u})\ldots{\phi}^{\lambda_\eta}_{\imath_\eta}(\mathfrak{u}),
\quad{\phi}^\emptyset_\imath\equiv1
\end{equation}
of variable  $u\in\mathfrak{U}$ belonging to $L^\infty_\chi$.
The system of $\chi$-essentially bounded  functions
\begin{equation*}\label{base2}
\phi^\mathbb{Y}:=\bigcup\big\{{\phi}^{\mathbb{Y}_n}\colon{n}\in\mathbb{N}_0\big\}
\quad\text{with}\quad
{\phi}^{\mathbb{Y}_n}:=\bigcup\big\{{\phi}^\lambda_\imath\colon
\imath^\lambda\vdash n, \imath\in\mathscr{I}_\eta\big\}
\end{equation*}
 is orthogonal in the space $L^2_\chi$ and is normed as follows
\[
\|{\phi}^\lambda_\imath\|_\chi^2=
\int |{\phi}^\lambda_\imath|^2d\chi=\lambda!\beta_\lambda,
\quad\imath^\lambda\vdash n,\quad
\lambda!:=\lambda_1!\ldots\lambda_\eta!.
\]
\end{lemma}

\begin{proof}
According to \eqref {projective}, we have
$(\pi_m\circ\pi_{m+l}^{-1}) u_{m+l}(\mathfrak{e}_m)= u_m(\mathfrak{e}_m)$ for $t=-1$ and
$(\pi_m\circ\pi_{m+l}^{-1})u_{m+l}(\mathfrak{e}_m)=u_m(\mathfrak{e}_m)-[a(1+t)^{-1}b]\mathfrak{e}_m$ for $t\noteq -1$ for any integer $l\ge1$. This means that
${({\phi}_k\circ\pi_m^{-1})(u_m)}=
\left\langle u_m(\mathfrak{e}_m) \mid \mathfrak{e}_k\right\rangle\nequiv0$ for all  $k\le m$ and that
\begin{equation}\label{stabl}
\begin{split}
{(\phi_m\circ\pi_{m+l}^{-1})(u_{m+l})}&=
\left\langle u_m(\mathfrak{e}_m) \mid \mathfrak{e}_m\right\rangle\quad\text{for}\quad  t=-1, \\
{(\phi_m\circ\pi_{m+l}^{-1})(u_{m+l})}&={\left\langle u_m(\mathfrak{e}_m) \mid \mathfrak{e}_m\right\rangle}-
{a(1+t)^{-1}b\left\langle\mathfrak{e}_m\mid \mathfrak{e}_m\right\rangle}\quad\text{for}\quad t\noteq -1.
\end{split}\end{equation}

Let $U(\eta)$ with $\eta=\eta(\lambda)$ be the unitary group acting over the linear complex
$\mathop{\rm span}\left\{\mathfrak{e}_{\imath_1},\ldots,\mathfrak{e}_{\imath_\eta}\right\}$ in $H$.
Let  $\chi_\eta$ be the probability Haar measure on $U(\eta)$ and $\pi_\eta\colon\mathfrak{U}\to U(\eta)$ be the corresponding
projector. Using \eqref{occurs0} and \eqref{stabl}, we obtain
\begin{equation}\label{weig}
\begin{split}
\int|{\phi}^\lambda_\imath(\mathfrak{u})|^2d\chi(\mathfrak{u})&
=\lim\int|({\phi}^\lambda_\imath\circ\pi_m^{-1})(u_m)|^2d\chi_m(u_m) \\
&=\lim\int|({\phi}^{\lambda_1}_{\imath_1}\circ\pi_m^{-1})(u_m)\ldots
({\phi}^{\lambda_\eta}_{\imath_\eta}\circ\pi_m^{-1})(u_m)|^2d\chi_m(u_m)\\
&=\int|({\phi}^{\lambda_1}_{\imath_1}\circ\pi_\eta^{-1})(u_\eta)\ldots
({\phi}^{\lambda_\eta}_{\imath_\eta}\circ\pi_\eta^{-1})(u_\eta)|^2d\chi_\eta(u_\eta).
\end{split}\end{equation}
By \eqref{weig} and the known integral formula for  unitary groups $U(\eta)$ \cite[1.4.9]{RudinFT80}, we get
 \[
\int|{\phi}^\lambda_\imath|^2d\chi
=\int\prod_{k=1}^{\eta(\lambda)}\left|\left\langle u_\eta(\mathfrak{e}_\eta)\mid \mathfrak{e}_{\imath_k}\right\rangle\right|^2d\chi_\eta(u_\eta)=
\frac{(\eta(\lambda)-1)!\lambda!}{(\eta(\lambda)-1+|\lambda|)!}.
\]

On the other hand, the invariant property \eqref{inv0} provides the formula
\begin{equation}\label{inv2}
\int f\,d\chi=\frac{1}{2\pi}\int\!d\chi(\mathfrak{u})
\int_{-\pi}^{\pi}f\left[\exp(\mathbbm{i}\vartheta)\mathfrak{u}\right]d\vartheta,
\qquad  f\in  L_\chi^\infty.
\end{equation}
From \eqref{inv2} it  follows the orthogonality relations
${\phi^{\lambda'}_\jmath\perp\phi^\lambda_\imath}$ with
${|\lambda'|\noteq |\lambda|}$, since
\[
\int{\phi}^{\lambda'}_\jmath\bar{{\phi}}^\lambda_\imath\,d\chi=
\frac{1}{2\pi}\int\phi^{\lambda'}_\jmath\bar{{\phi}}^\lambda_\imath\,d\chi
\int_{-\pi}^\pi{\exp\left[\mathbbm{i}\big(|\lambda'|-|\lambda|\big)\vartheta\right]}\,d\vartheta=0
\]
for any  $\lambda',\lambda\in\mathbb{Y}\setminus\{\emptyset\}$. Let $|\lambda'|=|\lambda|$ and
$\eta(\lambda')>\eta(\lambda)$
for definiteness. Then there exists an index $k$ with a nonzero integer  $\lambda'_k$ in
$\lambda'=\big(\lambda'_1,\ldots,\lambda'_k,\ldots,\lambda'_{\eta(\lambda')}\big)\in\mathbb{Y}\setminus\{\emptyset\}$
such that $\eta(\lambda)<k\le\eta(\lambda')$. In this case
${{\phi}^{\lambda'}_\jmath\perp{\phi}^\lambda_\imath}$
because \eqref{inv2} yields
\[
\int{\phi}^{\lambda'}_\jmath\bar{{\phi}^\lambda_\imath}\,d\chi=
\frac{1}{2\pi}\int{\phi}^{\lambda'}_\jmath\bar{{\phi}}^\lambda_\imath\,d\chi
\int_{-\pi}^{\pi}\exp\left(\mathbbm{i}\lambda'_k\vartheta\right)d\vartheta=0.
\]
Consider the case $|\lambda'|=|\lambda|$ and $\eta(\lambda')=\eta(\lambda)$.
If ${{\phi}^{\lambda'}_\jmath\noteq{\phi}^\lambda_\imath}$ then $\lambda'\noteq\lambda$.
There exists an index $0<k\le\eta(\lambda)$  such that $\lambda'_k\noteq\lambda_k$. As above,
${{\phi}^{\lambda'}_\jmath\perp{\phi}^\lambda_\imath}$, because
\[
\int{\phi}^{\lambda'}_\jmath\bar{{\phi}}^\lambda_\imath\,d\chi=
\frac{1}{2\pi}\int{\phi}^{\lambda'}_\jmath\bar{{\phi}}^\lambda_\imath\,d\chi
\int_{-\pi}^\pi\exp\left[\mathbbm{i}(\lambda'_k-\lambda_k)\vartheta\right]d\vartheta=0.
\]
This proves that the system ${\phi}^\mathbb{Y}$ is orthogonal.
\end{proof}

\section{Orthonormal basis of Schur polynomials}\label{4}

Let  $\imath^\lambda\vdash n$, $\eta=\eta(\lambda)$ and
$t_\imath=(t_{\imath_1},\ldots,t_{\imath_\eta})$ be a complex variable. Let
$t^\lambda_\imath:=\prod t_{\imath_j}^{\lambda_j}$.
The $n$-homogenous  Schur polynomial  is defined (see, e.g. \cite{MacDonald98}) to be
\[
s^\lambda_\imath(t_\imath):={D_\lambda(t_\imath)}/{\Delta(t_\imath)}\quad\text{where}\quad
D_\lambda(t_\imath)=\mathop{\mathrm{det}}\big[t_{\imath_i}^{\lambda_j+\eta-j}\big]
\text{ with }\lambda_j=0 \text{ for }j>\eta,
\]
$\Delta(t_\imath)={\prod_{1\le i<j\le\eta}(t_{\imath_i}-t_{\imath_j})}$ is
Vandermonde's determinant. It can be written as  $s^\lambda_\imath(t_\imath)=
{\sum}_{[\imath^\lambda]}t^\lambda_\imath$ with
summation over all semistandard Young tabloids \cite[I.2.2]{Fulton97}.

We construct an orthonormal basis in $L^2_\chi$ consisting of Schur polynomials on Paley--Wiener maps.
Assign (uniquely) to $\imath\in\mathscr{I}_\eta$ the vector ${\phi}_\imath:=\big({\phi}_{\imath_1},\ldots,{\phi}_{\imath_\eta}\big)$.
Let $s^\lambda_\imath(\mathfrak{u})=(s^\lambda_\imath\circ{\phi}_\imath)(\mathfrak{u})$
be $n$-homogeneous  functions of variable ${\mathfrak{u}\in\mathfrak{U}}$
with $\lambda\in\mathbb{N}^\eta$, defined by the formulas \eqref{schur0}.
Denote
\[
s^\mathbb{Y}_n:=\bigcup\big\{s^\lambda_\imath\colon
\imath^\lambda\vdash n\big\}, \quad
s^\mathbb{Y}:=\bigcup\big\{s^\mathbb{Y}_n\colon{n}\in\mathbb{N}_0\big\}\quad\text{with}
\quad{s_0=s^\emptyset_\imath\equiv1}.
\]

\begin{theorem}\label{Schur1}
The system of Schur polynomials $s^\mathbb{Y}$ forms an orthonormal
basis in $L^2_\chi$ and $s^\mathbb{Y}_n$ is the same basis in $L_\chi^{2,n}$.
The following  orthogonal decomposition holds,
\begin{equation}\label{ort}
L^2_\chi=\mathbb{C}\oplus L_\chi^{2,1}\oplus L_\chi^{2,2}\oplus\ldots.
\end{equation}
For any ${h\in{H}}$ the equality
 \eqref{trace} uniquely defines the conjugate-linear embedding
\begin{equation}\label{P_W}
{\phi}\colon{H}\ni{h}\longmapsto{\phi}_h\in L^2_\chi\quad\text{such that}\quad
\|{\phi}_h\|_\chi=\|h\|.
\end{equation}
\end{theorem}

\begin{proof}
Let $U(\eta)$ be the unitary group over the linear complex
$\mathop{\rm span}\left\{\mathfrak{e}_{\imath_1},\ldots,\mathfrak{e}_{\imath_\eta}\right\}$
with $\eta=\eta(\lambda)$. Taking into account \eqref{stabl} similarly as \eqref{weig}, we obtain
\[
\int{s}_\imath^\lambda\bar{s}_\imath^\mu\,d\chi=\int{s}_\imath^\lambda(z_\eta)
\,\bar{s}_\imath^\mu(z_\eta)\,d\chi_\eta(z_\eta)=\delta_{\lambda\mu}
\]
for all $[\imath^\lambda]$, $[\imath^\mu]$ with $\imath=(\imath_1,\ldots,\imath_\eta)$ and $\lambda,\mu\in\mathbb{N}^\eta$. In fact,
the corresponding Schur polynomials $\big\{s^\lambda_\imath\colon\lambda\in\mathbb{N}^\eta\big\}$  are  characters of the group
$U(\eta)$. Hence, by the Weyl integration formula, the right-hand side  integral is equal to Kronecker's delta $\delta_{\lambda\mu}$
\cite[Thm 8.3.2 \& Thm 11.9.1]{procesi2007lie}.

The family of finite alphabets ${\imath\in\mathscr{I}}$ is directed and for any $\imath,\imath'$
there exists $\imath''$ such that $\imath\cup\imath'\subset\imath''$. This means that the whole system $s^\mathbb{Y}_n$ is orthonormal in $L^2_\chi$.

The property ${s_\jmath^\mu\perp s_\imath^\lambda}$ with
${|\mu|\noteq |\lambda|}$ for any ${\imath,\jmath\in\mathscr{I}}$  follows from \eqref{inv2}, since
\[
\int{s}_\jmath^\mu\bar s_\imath^\lambda\,d\chi=\frac{1}{2\pi}\int{s}_\jmath^\mu\bar s_\imath^\lambda\,d\chi
\int_{-\pi}^\pi{\exp\big(\mathbbm{i}(|\mu|-|\lambda|)\vartheta\big)}\,d\vartheta=0
\]
for all  $\lambda\in\mathbb{Y}$ and $\mu\in\mathbb{Y}\setminus\{\emptyset\}$.
This yields  $L_\chi^{2,|\mu|}\perp L_\chi^{2,|\lambda|}$ in the space $L^2_\chi$.
Taking $\lambda=\emptyset$ with $|\emptyset|=0$, we get $1\perp L_\chi^{2,|\mu|}$ for all
$\mu\in\mathbb{Y}\setminus\{\emptyset\}$. Hence, \eqref{ort} is proved.

By Lemma~\ref{proh2} the subsystem  $\phi_k=s^1_k$ is orthonormal in $L^2_\chi$, hence
by Definition~\ref{PWm} it instantly follows that
$\|{\phi}_h\|_\chi^2=\sum|\mathfrak{e}_k^*(h)|^2\int|\phi_k|^2d\chi=\|h\|^2.$ It follows the isometric embedding \eqref{P_W}.

The set $\check{U}(m)$ of matrices with no eigenvalue $\{-1\}$ has Stone--\^{C}ech compactification  $\tilde{U}(m)$ such that the mapping $\check{\pi}^{m+1}_m$ has a continuous ${U}(m)$-valued extension  \[\tilde{\pi}^{m+1}_m\colon\tilde{U}(m+1)\longrightarrow{U}(m).\]
This fact follows from \cite[Thm 19.5]{Willard04} by virtue of that ${U}(m)$ is compact.
Hence, the projective limit $\tilde{\mathfrak{U}}:=\varprojlim\tilde{U}(m),$
determined by $\tilde{\pi}^{m+1}_m$,  is a compact set in $\mathfrak{U}$ with continuous ${U}(m)$-valued projections $\tilde\pi_m\colon\tilde{\mathfrak{U}}\rightarrow{U}(m)$.

Since $U(\infty)$  on $H$ acts irreducibly, for any $\mathfrak{u}'\noteq\mathfrak{u}''$
there is $m$ such that
\[
\phi_m(\mathfrak{u}')=\left\langle\pi_m(\mathfrak{u}')(\mathfrak{e}_m)\mid\mathfrak{e}_m\right\rangle\noteq
\left\langle\pi_m(\mathfrak{u}'')(\mathfrak{e}_m)\mid\mathfrak{e}_m\right\rangle=\phi_m(\mathfrak{u}''),
\]
i.e., $\phi^\mathbb{Y}$ separates $\mathfrak{U}$ and so $\tilde{\mathfrak{U}}$.
Hence, the system of Schur polynomials $s^\mathbb{Y}$ also separates
$\tilde{\mathfrak{U}}$. Moreover, each complex-conjugate function $\bar\phi_m(\mathfrak{u})=
\left\langle\mathfrak{e}_m\mid\pi_m(\mathfrak{u})(\mathfrak{e}_m)\right\rangle=
\left\langle\pi_m(\mathfrak{u}^\star)(\mathfrak{e}_m)\mid\mathfrak{e}_m\right\rangle$ belongs to
$\phi^\mathbb{Y}$. Thus, by the Stone--Weierstrass approximation theorem
the complex linear span of polynomials  $\phi^\mathbb{Y}$, as well as, of $s^\mathbb{Y}$,
forms a dense subspace in the Banach space of all continuous functions $C(\tilde{\mathfrak{U}})$.

Let $\tilde\chi_m$ means the image of $\chi_m$ under $\check{U}(m)\looparrowright{U}(m)$.
In Lemma~\ref{proh0} it inductively was shown that for every $\varepsilon>0$ there exists a compact set
$\varprojlim K_m\subset\check{\mathfrak{U}}$ such that
\[\tilde{\chi}_m(K_m)\ge1-\varepsilon \quad\text{for all}\quad m\]
where
$\tilde{\chi}_m(K_m)=\check{\chi}_m(K_m)=\chi_m(K_m)$, by definition of the measure  $\tilde\chi_m$ as an image. Hence, by the Prokhorov theorem
the projective limit $\tilde{\chi}=\varprojlim\tilde{\chi}_m$, defined by mappings $\tilde{\pi}^{m+1}_m$, possesses the properties
\[\tilde\chi(\Omega)=\inf\tilde{\chi}_m(\Omega)=\inf\chi_m(\Omega)=\varprojlim\chi_m(\Omega)=
\chi(\Omega)\] for all Borel $\Omega$ in $\check{\mathfrak{U}}$ or otherwise
$\tilde\chi|_{\check{\mathfrak{U}}}=\chi|_{\check{\mathfrak{U}}}$. Consequently,
\[\tilde\chi|_{\check{\mathfrak{U}}}=
\chi|_{\check{\mathfrak{U}}}=\chi|_{\check{\mathfrak{U}}\bigsqcup(\mathfrak{U}\setminus\check{\mathfrak{U}})}
=\chi|_\mathfrak{U}\quad \text{since}\quad \chi(\mathfrak{U}\setminus\check{\mathfrak{U}})=0.\]
In particular, $\tilde{\chi}=\varprojlim\tilde{\chi}_m$ is regular on $\tilde{\mathfrak{U}}$  by
the Riesz--Markov theorem \cite[1.1]{Nelson64}.

As a consequence, the space $L^2_\chi$ coincides with the completion of $C(\tilde{\mathfrak{U}})$ and
for any ${f\in L_\chi^2}$  there exists a sequence ${(f_n)\subset\mathop{\rm span}(s^\mathbb{Y})}$ such that
${\int|f-f_n|^2d\chi\to0}$. Hence, the system $s^\mathbb{Y}$ forms an orthogonal basis in $L_\chi^2$.

Finally,  $s^\mathbb{Y}_n\cap L^2_\chi$ is total in  $L_\chi^{2,n}$
and $s^\mathbb{Y}_n\perp s^\mathbb{Y}_m$ if $n\noteq m$. This yields \eqref{ort}.
\end{proof}

\section{Unitarily-weighted symmetric Fock space}\label{5}

Define on the tensor power ${H}^{\otimes n}$ the unitarily-weighted norm
$\|\cdot\|_{H^{\otimes  n}_\beta}=
{\langle\cdot\mid\cdot\rangle^{1/2}_{H^{\otimes  n}_\beta}}$ where the
inner product ${\langle\cdot\mid\cdot\rangle^{1/2}_{H^{\otimes  n}_\beta}}$ is determined by the relations
\begin{equation}\label{deffull}
\langle\mathfrak{e}^{\otimes\lambda}_\imath\mid
\mathfrak{e}^{\otimes\lambda'}_{\imath'}\rangle_{H^{\otimes  n}_\beta}
=\left\{
\begin{array}{clcl}
\dfrac{(\eta-1)!}{(\eta-1+n)!}&: \lambda=\lambda'&\text{and}&\imath=\imath' \\
  0 &: \lambda\noteq\lambda'&\text{or}&\imath\noteq\imath'.
\end{array}\right.
\end{equation}
Here $\mathfrak{e}^{\otimes\lambda}_\imath:=\sigma'(\mathfrak{e}^{\otimes\lambda_1}_{\imath_1}\otimes\ldots\otimes
\mathfrak{e}^{\otimes\lambda_\eta}_{\imath_\eta})$  with $\eta=\eta(\lambda)$ and $\sigma'\in S_n$ is fixed.
Let ${H}^{\otimes n}_\beta$ be the completion of $\big\{{H}^{\otimes n},\|\cdot\|_{H^{\otimes  n}_\beta}\big\}$.
Its closed subspace, defined by the projection
\[{\mathcal{S}_n\colon{H}^{\otimes n}_\beta\ni\mathfrak{e}^{\otimes\lambda}_\imath
\longmapsto\mathfrak{e}^{\odot\lambda}_\imath=(n!)^{-1}{\sum}_{\sigma\in S_n}\sigma(\mathfrak{e}^{\otimes\lambda}_\imath)}\]
forms an unitarily-weighted symmetric tensor power ${H}^{\odot n}_\beta\subset{H}^{\otimes n}_\beta$
with the inner product determined by relations
$\langle\mathfrak{e}^{\odot\lambda}_\imath\mid
\mathfrak{e}^{\odot\lambda'}_{\imath'}\rangle_{H^{\otimes  n}_\beta}=
\beta_\lambda\langle\mathfrak{e}^{\odot\lambda}_\imath\mid
\mathfrak{e}^{\odot\lambda'}_{\imath'}\rangle$ or more specific
\begin{equation}\label{defw}
\langle\mathfrak{e}^{\odot\lambda}_\imath\mid
\mathfrak{e}^{\odot\lambda'}_{\imath'}\rangle_{H^{\otimes  n}_\beta}
=\left\{
\begin{array}{clcl}
\dfrac{\lambda!}{n!}\dfrac{(\eta-1)!}{(\eta-1+n)!}&: \lambda=\lambda&\text{and}&\imath=\imath' \\
  0 &: \lambda\noteq\lambda'&\text{or}&\imath\noteq\imath'.
\end{array}\right.
\end{equation}

\begin{definition}
The {\em unitarily-weighted symmetric Fock space} is defined to be the Hilbertian orthogonal sum
 ${\Gamma_{\!\beta}(H)=\bigoplus_{n\ge0}{H}_\beta^{\odot n}}$ of elements $\psi=\bigoplus\psi_n$, $\psi_n\in{H}_\beta^{\odot n}$
with the orthogonal basis $\mathfrak{e}^{\mathbb{Y}}=
\bigcup\big\{\mathfrak{e}^{\mathbb{Y}_n}\colon{n}\in\mathbb{N}_0\big\}$ and the following inner product and norm
\[
\langle\psi\mid\psi'\rangle_\beta
=\sum{n!}\langle\psi_n\mid\psi_n'\rangle_{H^{\otimes n}_\beta},\quad\|\psi\|_\beta=\langle\psi\mid\psi\rangle_\beta^{1/2}.
\]
\end{definition}

We immediately notice that $\|h\|_\beta^2=\sum|\mathfrak{e}_i^*(h)|^2=\|h\|^2$ for all $h=\sum\mathfrak{e}_i\mathfrak{e}_i^*(h)\in{H}$.

\begin{lemma}\label{T_a}
The set of coherent states $\left\{\exp(h)\colon{h}\in{H}\right\}$ is total in
 $\Gamma_\beta(H)$ and the expansion \eqref{Tayl} is convergent in  $\Gamma_\beta(H)$.
The injections \[\Gamma(H)\looparrowright\Gamma_\beta(H)\quad\text{and}\quad {{H}^{\odot n}
\looparrowright{H}^{\odot n}_\beta}\] are contractive and dense. The
$\Gamma_\beta(H)$-valued function ${H\ni{h}\longmapsto\exp(h)}$ is entire analytic.
The shift group, defined to be
\[
\mathcal{T}_a\exp(h):=\exp(h+a)
=\exp(\partial_a)\exp(h)\quad\text{with}\quad
\partial_a\exp(h)=\frac{d\exp(h+za)}{dz}{\Big|}_{z=0}
\]
for $a,h\in H$, has a unique linear extension ${\mathcal{T}_a\colon
\Gamma_\beta(H)\ni\psi\longmapsto\mathcal{T}_a\psi\in\Gamma_\beta(H)}$ such that
 \begin{equation}\label{contract}
 \|\mathcal{T}_{a}\psi\|^2_\beta\le\exp\big(\|a\|^2\big)\|\psi\|^2_\beta
 \quad\text{and}\quad\mathcal{T}_{a+b}=
\mathcal{T}_a\mathcal{T}_b=\mathcal{T}_b\mathcal{T}_a,\quad{a,b\in H}.
\end{equation}
\end{lemma}

\begin{proof}
Taking  into account that $\beta_\lambda\le1$, we get the following inequalities
\begin{align*}
\|h^{\otimes n}\|^2_{H^{\otimes  n}_\beta}&\!=\!\sum_{\imath^\lambda\vdash n}\!
\Big(\frac{n!}{\lambda !}\Big)^2\|\mathfrak{e}^{\odot\lambda}_\imath\|^2_{H^{\otimes n}_\beta}
|\mathfrak{e}^{*\lambda}_\imath(h)|^2\!=\!\sum_{\imath^\lambda\vdash n}\beta_\lambda
\frac{n!}{\lambda !}|\mathfrak{e}^{*\lambda}_\imath(h)|^2\!
\le\|h^{\otimes n}\|^2 =\|h\|^{2n},\\
&\|\exp(h)\|^2_\beta
=\sum_{n\ge0}\frac{1}{n!}\sum_{\imath^\lambda\vdash n}\beta_\lambda\frac{n!}{\lambda !}
|\mathfrak{e}^{*\lambda}_\imath(h)|^2
\stackrel{\eqref{constant}}\le\exp\|h\|^2\stackrel{\eqref{exp}}=\|\exp(h)\|^2_\Gamma.
\end{align*}
Hence, \eqref{Tayl0},  \eqref{Tayl} are convergent in $\Gamma_\beta(H)$. This implies that ${h\mapsto\exp(h)}$
is analytic and inclusions
$\Gamma(H)\looparrowright\Gamma_\beta(H)$ and
${{H}^{\odot n}\looparrowright{H}_\beta^{\odot n}}$ are contractive. By the polarization formula \eqref{polarization} their  ranges are dense.

Using the binomial formula  ${(h+za)^{\otimes n}}=
\bigoplus_{m=0}^{n}\binom{n}{m}(za)^{\otimes m}\odot h^{\otimes (n-m)}$, we find
\[
\partial_{a}^m\exp(h)=
\frac{d^m\exp(h+za)}{dz^m}{\mathrel{\Big|}}_{z=0}
=\bigoplus_{n\ge m}\frac{\mathcal{S}_{n/m}[a^{\otimes m}\otimes h^{\otimes (n-m)}]}{(n-m)!},\quad
{z\in\mathbb{C}}
\]
with the orthogonal projector $\mathcal{S}_{n/m}$ defined as
$\psi_m\odot\psi_{n-m}=\mathcal{S}_{n/m}\left(\psi_m\otimes\psi_{n-m}\right)\in{H}^{\odot n}_\beta$
for all $\psi_m\in{H}^{\odot m}_\beta$ and $\psi_{n-m}\in{H}^{\odot(n- m)}_\beta$. By orthogonality
$\|\mathcal{S}_{n/m}\|\le1$.

Applying the expansions \eqref{Tayl0} to $a^{\otimes m}$ and $h^{\otimes(n-m)}$, by  \eqref{deffull}, we get
\[
\|a^{\otimes m}\otimes h^{\otimes(n-m)} \|^2_{H^{\otimes n}_\beta}=\!\!\sum_{ \substack{\imath^\lambda\vdash m \\  \jmath^\mu\vdash(n-m)}}\!\!
\Big(\frac{m!}{\lambda !}\frac{(n-m)!}{\mu!}\Big)^2
\|\mathfrak{e}^{\odot\lambda}_\imath\otimes\mathfrak{e}^{\odot\mu}_\jmath\|^2_{H^{\otimes  n}_\beta}
|\mathfrak{e}^{*\lambda}_\imath(a)|^2|\mathfrak{e}^{*\mu}_\jmath(h)|^2
\]
with summations over semistandard  tableaux $[\imath^\lambda],[\jmath^\mu]$  and $\imath,\jmath\in\mathscr{I}$.
Let $(\lambda,\mu)\in\mathbb{N}^{\eta(\lambda,\mu)}$ be the smallest partition
of number $n$ with the length $\eta(\lambda,\mu)$  containing the partitions $\lambda$ for $m$ and $\mu$ for $n-m$.
Then $\eta(\lambda,\mu)\ge\max\{\eta(\lambda),\eta(\mu)\}$ and so
\[
\|\mathfrak{e}^{\odot\lambda}_\imath\otimes\mathfrak{e}^{\odot\mu}_\jmath\|^2
_{H^{\otimes  n}_\beta}=\frac{(\eta(\lambda,\mu)-1)!}{(\eta(\lambda,\mu)-1+n)!}\le\min\{\beta_\lambda,\beta_\mu\},
\]
since $\frac{(\eta-1)!}{(\eta-1+n)!}$ is decreasing in variable $\eta$. Thus, the following inequality
\begin{align*}
\|a^{\otimes m}\otimes h^{\otimes(n-m)} \|^2_{H^{\otimes n}_\beta}&\le
\sum_{ \substack{\imath^\lambda\vdash m \\  \jmath^\mu\vdash(n-m)}}
\Big(\frac{m!}{\lambda !}\frac{(n-m)!}{\mu!}\Big)^2\min\{\beta_\lambda,\beta_\mu\}
|\mathfrak{e}^{*\lambda}_\imath(a)|^2|\mathfrak{e}^{*\mu}_\jmath(h)|^2\\
&=\|a^{\otimes m}\|^2\|h^{\otimes(n-m)} \|^2_{H^{\otimes(n-m)}_\beta}=\|a\|^{2m}\|h^{\otimes(n-m)} \|^2_{H^{\otimes(n-m)}_\beta}
\end{align*}
holds. Using this inequality and  that $\|\mathcal{S}_{n/m}\|\le1$, we find
\[\begin{split}
\|\partial_{a}^m\exp(h)\|^2_\beta&=\!\sum_{n\ge m}\!\!
\frac{\|\mathcal{S}_{n/m}[a^{\otimes m}\otimes h^{\otimes (n-m)}]\|^2_\beta}{(n-m)!}
\le\!\sum_{n\ge m}\!\!\frac{\|\mathcal{S}_{n/m}\|^2\|a^{\otimes m}\otimes h^{\otimes (n-m)}\|^2_\beta}{(n-m)!}\\
&\le\|{a}^{\otimes m}\|^2\sum_{n\ge m}
\frac{\|\mathcal{S}_{n/m}\|^2\|{h}^{\otimes (n-m)}\|^2_\beta}{(n-m)!}
\le\|{a}\|^{2m}\|\exp(h)\|^2_\beta.
  \end{split}\]
Summing with coefficients $1/m!$, we get
$\|\mathcal{T}_{a}\exp(h)\|^2_\beta\le\exp\big(\|a\|^2\big)\|\exp(h)\|^2_\beta$.
This inequality and totality of $\left\{\exp(x)\colon h\in{H}\right\}$ in
$\Gamma_\beta(H)$ yield the required inequality \eqref{contract}. It also follows that
$\Gamma_\beta(H)$ is invariant under $\mathcal{T}_{a}$ and that the group property \eqref{contract} holds,
since $\partial_{a+b}=\partial_a+\partial_b$ for all $a,b\in H$ by linearity.
\end{proof}

\begin{lemma}\label{Paley--Wiener}
The mapping  $\phi\colon H\ni{h}\longmapsto \phi_h\in L^2_\chi $, extended onto
$\mathcal{T}_{a}\exp(h)$ as
\[
\varPhi\colon\mathcal{T}_{a}\exp(h)\longmapsto\sum_{n\ge0}\frac{1}{n!}\sum_{\imath^\lambda\vdash n}
\frac{n!}{\lambda!}\phi_\imath^\lambda\mathfrak{e}^{*\lambda}_\imath(h+a),\quad
a\in H,
\]
has the unique isometric conjugate-linear extension
\[{\varPhi\colon\Gamma_\beta(H)\ni\psi\longmapsto\varPhi\psi\in L_\chi^2 }\quad\text{
with the adjoint mapping}\quad {\varPhi^*\colon L_\chi^2 \rightarrow\Gamma_\beta(H)}\] defined to be
$\langle\varPhi\mathfrak{e}^{\odot\lambda}_\imath\mid f\rangle_\chi=
\langle\mathfrak{e}^{\odot\lambda}_\imath\mid\varPhi^*f\rangle_\beta$ for all
${f\in L_\chi^2 }$ in such way that
\[\varPhi\colon
\mathfrak{e}^{\odot \lambda}_\imath/\|\mathfrak{e}^{\odot \lambda}_\imath\|_\beta
\longmapsto\phi_\imath^\lambda/\|{\phi}^\lambda_\imath\|_\chi\quad\text{for all}\quad
 \lambda\in\mathbb{Y}, \ \imath\in\mathscr{I}_{\eta(\lambda)}.\]
 As a result, the conjugate-linear isometries $\Gamma_\beta(H) \stackrel{\varPhi}\simeq L^2_\chi$ and
${H}^{\odot n}_\beta\stackrel{\varPhi}\simeq L^{2,n}_\chi$ hold.
\end{lemma}

\begin{proof}
By Lemma~\ref{T_a}
the $\Gamma_\beta(H)$-valued function $H\ni h\mapsto\mathcal{T}_{a}\exp(h)$ is well defined
 for all ${a\in H}$. Let us use the expansion $\phi_{h+a}={\sum\mathfrak{e}^*_i(h+a)\phi_i}$.
By Lemma~\ref{proh2} and Theorem~\ref{Schur1},
$\phi\colon H\ni{h}\longmapsto \phi_h\in L_\chi^2 $ may be extended to $\varPhi$ in following way
\[\begin{split}
\varPhi\mathcal{T}_{a}\exp(h)
&=\sum_{n\ge0}\frac{1}{n!}\sum_{\imath^\lambda\vdash n}
\frac{n!}{\lambda!}\phi_\imath^\lambda\mathfrak{e}^{*\lambda}_\imath(h+a)
=\prod_{i\ge0}\sum_{n\ge0}\frac{\phi_i^n}{n!}\mathfrak{e}^{*n}_i(h+a)\\
&=\prod\exp\left(\phi_i\mathfrak{e}^*_i(h+a)\right)=
\exp\left(\phi_{h+a}\right)\quad\text{where}\\
\varPhi[(h+a)^{\odot n}]&=\phi_{h+a}^n=\sum\limits_{\imath^\lambda\vdash n}
\frac{n!}{\lambda!}\phi_\imath^\lambda\mathfrak{e}^{*\lambda}_\imath(h+a),\quad a\in H
\end{split}\]
is an orthogonal component of $\varPhi\mathcal{T}_{a}\exp(h)$ in $L^2_\chi $. It  follows that
\[\begin{split}
\|\exp(\phi_{h+a})\|_\chi^2&=\sum_{n\ge0}\frac{1}{n!^2}
\sum_{\imath^\lambda\vdash n}\|\phi_\imath^\lambda\|_\chi^2
\frac{n!^2}{\lambda!^2}|\mathfrak{e}^{*\lambda}_\imath(h+a)|^2\\
&=\sum_{n\ge0}\frac{1}{n!^2}\sum_{\imath^\lambda\vdash n}
\frac{n!^2}{\lambda!}\beta_\lambda|\mathfrak{e}^{*\lambda}_\imath(h+a)|^2\le\sum_{n\ge0}\frac{1}{n!}\sum_{\imath^\lambda\vdash n}
\frac{n!}{\lambda!}|\mathfrak{e}^{*\lambda}_\imath(h+a)|^2\\
&=\prod\exp|\mathfrak{e}^*_i(h+a)|^2=\exp\|h+a\|^2.
\end{split}\]
Hence, the composition $\mathfrak{U}\ni\mathfrak{u}\longmapsto[\varPhi\exp(h+a)](\mathfrak{u})$ is well defined in $L^2_\chi $.

Now, we consider the ordinary irreducible representation of permutation group $S_n$ on the Specht $\lambda$-module
$S^\lambda_\imath$ that is corresponded to the standard Young tableau $[\imath^\lambda]$. The
following known hook formula (see \cite[I.4.3]{Fulton97}) holds,
\begin{equation}\label{hook}
\hbar_\lambda:={n!}\Big({\prod}_{i\le \lambda_j}h(i,j)\Big)^{-1}\quad\text{where}\quad
\hbar_\lambda=\mathop{\rm{dim}}S^\lambda_\imath,
\end{equation}
with $h(i,j)\!=\!\#\big\{\Box_{i'j'}\in[\imath^\lambda]\colon i'\ge i,j'=j\big\}\!=\!\#\big\{\Box_{i'j'}\in[\imath^\lambda]\colon {i'=i, j'\ge j}\big\}$ independed of $\imath\in\mathscr{I}$.
Assign to $\imath\in\mathscr{I}_\eta$ the vectors
\[\left({\phi}_{\imath_1}(\mathfrak{u}){\mathfrak{e}}^*_{\imath_1}(h),\ldots,
{\phi}_{\imath_\eta}(\mathfrak{u}){\mathfrak{e}}^*_{\imath_\eta}(h)\right):=t_\imath(\mathfrak{u},h).\] Let
$s^\lambda_\imath(\mathfrak{u},h):=s^\lambda_\imath(t_\imath)$
with $t_\imath=t_\imath(\mathfrak{u},h)$ for all ${\mathfrak{u}\in\mathfrak{U}}$,
where polynomial terms are
${\phi}^\lambda_\imath(\mathfrak{u}){\mathfrak{e}}_\imath^{*\lambda}(h)=
{\phi}^{\lambda_1}_{\imath_1}(\mathfrak{u}){\mathfrak{e}}^{*\lambda_1}_{\imath_1}(h)
\ldots{\phi}^{\lambda_\eta}_{\imath_\eta}(\mathfrak{u})
{\mathfrak{e}}^{*\lambda_\eta}_{\imath_\eta}(h)$.
Applying  the Frobenius formula \cite[I.7]{MacDonald98} and taking into account
  \eqref{trace}, \eqref{schur0}, \eqref{hook}, we obtain
\[
\phi_h^n(\mathfrak{u})=
{\sum}_{\imath^\lambda\vdash n}\hbar_\lambda
s^\lambda_\imath(\mathfrak{u},h),\quad h\in H
\]
where $s^\lambda_\imath=0$ if  $\lambda^\intercal_1>l_\lambda$
and  the summation is over all standard tabloids. Hence,
$\big\{\phi_h^n\colon h\in H\big\}$ is total in $L_\chi^{2,n}$
by Theorem~\ref{Schur1}. In consequence,
$\left\{\exp(\phi_h)\colon h\in H\right\}$ is total in $L_\chi^2 $.
This yields surjectivity of $\varPhi$ and of all its restrictions to ${H}^{\odot n}_\beta$.
\end{proof}

\begin{corollary}\label{tot}
The sets $\big\{\phi_h^n\colon{h}\in{H}\big\}$ in $L^{2,n}_\chi$ and $\left\{\exp\phi_h\colon{h}\in{H}\right\}$ in $L^{2}_\chi$
are total.
\end{corollary}

\section{Fourier analysis on virtual unitary matrices}\label{6}

Consider the isometry
${H}^{*\odot n}_\beta\stackrel{\mathcal{P}}\simeq{P}_\beta^n(H)$
(see e.g., \cite[1.6]{Floret97}), where the space ${P}_\beta^n(H)$
of unitarily-weighted $n$-homogeneous Hilbert--Schmidt polynomials  of variable $h\in H$ is defined to be
a restriction to the diagonal in  ${H\times\ldots\times H}$ of the $n$-linear forms ${\mathcal{P}\circ\psi_n}$
endowed with the norm $\|\psi_n^*\|_{P_\beta^n}=\|\psi_n\|_{H^{\otimes n}_\beta}$ where
\[
\psi_n^*(h):=\langle h^{\otimes n}\mid\psi_n\rangle_{H^{\otimes  n}_\beta}\simeq
\left\langle ({h},\ldots, {h})\mid\mathcal{P}\circ\psi_n\right\rangle,
\quad{\psi_n\in H^{\odot n}_\beta}.
\]

Let ${H}^2_\beta=\sum_{n\ge0}{P}_\beta^n(H)$ be the direct sum  of functions $\psi^*(h)=\sum\psi_n^*(h)$ of variable ${h\in H}$ with summands
$\psi_n^*=\mathcal{P}\circ\psi_n\in {P}_\beta^n(H)$ where $\psi=\sum\psi_n\in\Gamma_\beta(H)$.
Since the set $\{\exp(h)\colon h\in H\}$ is total in $\Gamma_\beta(H)$, elements of ${H}^2_\beta$ can be written as
\[
{H}^2_\beta=\left\{\psi^*(h)=
\left\langle \exp(h)\mid\psi\right\rangle_\beta\colon \psi=\sum\psi_n\in\Gamma_\beta(H)\right\}.
\]
The analyticity of $H\ni{h}\mapsto\psi^*(h)$ is a result of the composition $\exp(\cdot)$ and $\psi^*(\cdot)$.
\begin{definition} Let ${H}^2_\beta$ be defined as a Hardy space of unitarily-weighted Hilbert--Schmidt analytic functions
$\psi^*(h)$ of variable ${h\in H}$ endowed with the inner product
 \[
 {\langle\psi^*(\cdot)\mid\varphi^*(\cdot)\rangle_{{H}^2_\beta}}:={\left\langle\varphi\mid\psi\right\rangle_\beta}\quad\text{where}\quad
 \|\psi^*\|_{{H}^2_\beta}^2=\langle\psi^*(\cdot)\mid\psi^*(\cdot)\rangle_{{H}^2_\beta}=
 \sum n!\|\psi_n^*\|^2_{P_\beta^n}.
\]
\end{definition}

The conjugate-linear surjective isometry from ${H}^2_\beta$ onto
$\Gamma_\beta(H)$ is realized by the conjugate-linear mapping
\[{*\colon\Gamma_\beta(H)\ni\psi\longmapsto\psi^*\in {H}^2_\beta},\quad \psi=\sum\psi_n.\]

On the other hand, the correspondence $\varPhi\colon\mathfrak{e}^{\odot \lambda}_\imath
\rightleftarrows\phi_\imath^\lambda$ with $\lambda\in\mathbb{Y}$ and $\imath\in\mathscr{I}_{\eta(\lambda)}$
allows us to determine the conjugate-linear isometry from $\Gamma_\beta(H)$ onto $L^2_\chi$.  As a result,
the mapping \[\varPsi\colon {H}^2_\beta\ni\mathfrak{e}^{*\lambda}_\imath/\|\mathfrak{e}^{\odot\lambda}_\imath\|_\beta
\longmapsto\phi_\imath^\lambda/\|{\phi}^\lambda_\imath\|_\chi\in L^2_\chi\] defines  the surjective isometry
\[\varPsi\colon{H}^2_\beta\longrightarrow L^2_\chi \quad\text{and its adjoint}\quad
{\varPsi^*\colon L^2_\chi \longrightarrow{H}^2_\beta}.\]

\begin{lemma}\label{pbasis}
The systems of Hilbert--Schmidt polynomials of variable ${h\in{H}}$,
\[
\mathfrak{e}^{*\mathbb{Y}_n}:=\bigcup\big\{\mathfrak{e}^{*\lambda}_\imath\colon
\imath^\lambda\vdash n, \imath\in\mathscr{I}\big\} \quad\text{and}\quad
\mathfrak{e}^{*\mathbb{Y}}:=\bigcup\big\{\mathfrak{e}^{*\mathbb{Y}_n}\colon{n}\in\mathbb{N}_0\big\}
\]
where ${\mathfrak{e}^{*\emptyset}_\imath=1}$,
form orthogonal bases in ${P}_\beta^n(H)$ and ${H}^2_\beta$, respectively, such that
\[
\|\mathfrak{e}^{*\lambda}_\imath\|_{P_\beta^n}^2=\beta_\lambda
\|\mathfrak{e}^{\odot\lambda}_\imath\|^2=
\dfrac{(\eta(\lambda)-1)!}{(\eta(\lambda)-1+n)!}\frac{\lambda!}{n!},
\quad\imath^\lambda\vdash n.
\]
Every function ${\psi^*\in{H}^2_\beta}$ with $\psi\in\Gamma_\beta(H)$
has the expansion with respect to $\mathfrak{e}^{*\mathbb{Y}}$
\begin{equation}\label{hhtaylor}
\psi^*(h)=\left\langle
\exp(h)\mid\psi\right\rangle_\beta
=\sum_{n\ge0}\frac{1}{n!}\sum_{\imath^\lambda\vdash n}\frac{n!}{\lambda!}
\mathfrak{e}^{*\lambda}_\imath(h)
\big\langle\mathfrak{e}^{\odot\lambda}_\imath\mid\psi_n\big\rangle_\beta
\end{equation}
with summation in the inner sum over all semistandard tabloids ${[\imath^\lambda]}$
such that ${\imath^\lambda\vdash n}$.
Each  function ${\psi^*\in{H}^2_\beta}$ is entire  Hilbert--Schmidt analytic and can be also written as
\begin{align}\begin{split}\label{cformula}
&\psi^*(h)=\big\langle\psi^*(\cdot)\mid\exp\langle\cdot\mid h\rangle\big\rangle_{H^2_\beta}
=\big\langle\psi^*(\cdot)\mid{E}(\cdot,h)\big\rangle_{{H}^2_\beta},\quad
{\psi\in\Gamma_\beta(H)}\\
\text{where}\quad&
E(h',h):=|\exp\langle h'\mid h\rangle|^2\!/\exp\langle h\mid h\rangle\quad\text{for all}\quad  h\in{H}.
\end{split}\end{align}
The following  linear isometries, defined by linearization via coherent states,  hold
\begin{equation}\label{PSI}
{H}^2_\beta\stackrel{\varPsi}\simeq{L}^2_\chi ,\quad
{P}_\beta^n(H)\stackrel{\varPsi}\simeq {L}^{2,n}_\chi .
\end{equation}
\end{lemma}

\begin{proof}
Taking into account \eqref{Tayl} and \eqref{defw},  we conclude that every
$\psi^*\in{H}^2_\beta$ such that $\psi=\bigoplus\psi_n\in\Gamma_\beta(H)$ with
${\psi_n\in{H}^{\odot n}_\beta}$ has the following expansion
\begin{equation*}\label{wort}
\psi^*(h)
=\sum_{n\ge0}\frac{1}{n!}\sum_{\imath^\lambda\vdash n}\frac{n!}{\lambda!}\mathfrak{e}^{*\lambda}_\imath(h)
\langle\mathfrak{e}^{\odot\lambda}_\imath\mid\psi_n\rangle_\beta  \quad \text{where} \quad
\psi=\bigoplus_{n\ge0}\sum_{\imath^\lambda\vdash n}
\frac{\langle\mathfrak{e}^{\odot\lambda}_\imath\mid\psi_n\rangle_\beta}
{\|\mathfrak{e}^{\odot \lambda}_\imath\|^2_\beta}\mathfrak{e}^{\odot\lambda}_\imath.
\end{equation*}
On the other hand, in relative to the inner product $\langle\cdot\mid\cdot\rangle_\Gamma$,
we have
\[
\exp\langle h'\mid h\rangle=\bigoplus_{n\ge0}\frac{1}{n!}\sum_{\imath^\lambda\vdash n}
\frac{n!}{\lambda !}\mathfrak{e}^{*\lambda}_\imath(h')\,\bar{\mathfrak{e}}^{*\lambda}_\imath(h)=
\sum_{n\ge0}\frac{1}{n!}\sum_{\imath^\lambda\vdash n}
\frac{\mathfrak{e}^{*\lambda}_\imath(h')\bar{\mathfrak{e}}^{*\lambda}_\imath(h)}
{\|\mathfrak{e}^{\odot \lambda}_\imath\|^2}.
\]
Verify the first equality  in \eqref{cformula} by substituting \eqref{hhtaylor} into the  formula \eqref{cformula}. We get
\begin{align*}
\psi^*(h)&=\bigg\langle\sum_{n\ge0}\sum_{\imath^\lambda\vdash n}
\frac{\langle\mathfrak{e}^{\odot\lambda}_\imath\mid\psi_n\rangle_\beta}
{\|\mathfrak{e}^{\odot \lambda}_\imath\|^2_\beta}\mathfrak{e}^{*\lambda}_\imath(h')\mid
\sum_{n\ge0}\frac{1}{n!}\sum_{\imath^\lambda\vdash n}
\frac{\mathfrak{e}^{*\lambda}_\imath(h')\bar{\mathfrak{e}}^{*\lambda}_\imath(h)}
{\|\mathfrak{e}^{\odot \lambda}_\imath\|^2}\bigg\rangle_{\!\!H^2_\beta}\\
&=\sum_{n\ge0}\frac{1}{n!}\sum_{\imath^\lambda\vdash n}\frac{n!}{\lambda!}
\mathfrak{e}^{*\lambda}_\imath(h)
\langle\mathfrak{e}^{\odot\lambda}_\imath\mid\psi_n\rangle_\beta
=\left\langle\exp(h)\mid\psi\right\rangle_\beta.
\end{align*}
If ${\omega^*(h'):=\psi^*(h)\exp\langle h\mid h'\rangle[\exp\langle h'\mid h'\rangle ]^{-1}}$
then $\omega^*(h)=\psi^*(h)$ for  ${h=h'\in{H}}$. Now, putting  $\omega^*(h'):=
\big\langle\psi^*(\cdot)\mid \exp\langle h'\mid\cdot\rangle [\exp\langle h'\mid h'\rangle ]^{-1}
\exp\langle \cdot\mid h'\rangle \big\rangle_{H^2_\beta}$, we obtain
\[\begin{split}
\psi^*(h)&=\omega^*(h)=\left\langle\omega^*\mid \exp(\cdot\mid h)
\right\rangle_{H^2_\beta}\\
&=\big\langle\psi^*(\cdot)\mid\exp(h\mid\cdot)[\exp(h\mid h)]^{-1}\exp(\cdot\mid h)
\big\rangle_{H^2_\beta}
=\big\langle\psi^*(\cdot)\mid{E}(\cdot,h)\big\rangle_{H^2_\beta}.
\end{split}\]
Hence, the second equality in \eqref{cformula} holds. Lemma~\ref{Paley--Wiener}  yields
\eqref{PSI}.
\end{proof}

\begin{remark}
Since  $\phi_h=\sum\mathfrak{e}^*_i(h)\phi_i$ for all $h=\sum\mathfrak{e}^*_i(h)\mathfrak{e}_i$,  a range of the embedding \eqref{P_W} coincides with $L_\chi^{2,1}$.
\end{remark}

\begin{lemma}\label{infty2}
Denote $\exp\langle h'\mid h\rangle:=K(h',h)$. The functions
\[H\ni{h}\longmapsto(\varPsi\circ{K})(\mathfrak{u},h)\quad  \text{and}\quad H\ni{h}\longmapsto(\varPsi\circ{E})(\mathfrak{u},h)\]
with $\mathfrak{u}\in\mathfrak{U}$
take values in $L^2_\chi $ and can be represented as follows
\[
(\varPsi\circ{K})(\mathfrak{u},h)=\exp\left(\phi_h(\mathfrak{u})\right),\qquad
(\varPsi\circ{E})(\mathfrak{u},h)=
\exp\big(2\mathop{\mathrm{Re}}\phi_h(\mathfrak{u})-\|h\|^2\big)
\]
where the last exponential function has the power series expansion
\begin{equation}\label{cp}
\begin{split}
\exp\left\{2\mathop{\mathrm{Re}}\phi_h-\|h\|^2\right\}&=
\sum_{m,n\ge0}\frac{\|h\|^{m+n}}{m!n!}
\mathfrak{h}_{n,m}\left(\phi_{h/\|h\|},\bar\phi_{h/\|h\|}\right)\\
\mathfrak{h}_{n,m}(z,\bar{z})&=\sum^{m\wedge n}_{k=0}(-1)^kk!\binom{m}{k}\binom{n}{k}{z}^{m-k}\bar{z}^{n-k}
\end{split}
\end{equation}
with coefficients in the form of complex Hermite polynomials $\mathfrak{h}_{n,m}(z,\bar{z})$, $z\in\mathbb{C}$.
\end{lemma}

\begin{proof}
Applying the transform $\varPsi$ to  ${K}(h',h)$ in variable ${h'\in{H}}$, we obtain
\begin{align*}
(\varPsi\circ{K})(\mathfrak{u},h)&=\!\sum_{n\ge0}\frac{1}{n!}
\sum\limits_{\imath^\lambda\vdash n}\frac{n!}{\lambda!}
\phi^\lambda_\imath(\mathfrak{u}){\mathfrak{e}}^{*\lambda}_\imath(h)
=\!\sum_{n\ge0}\frac{1}{n!}\Big(\sum_{i\ge0}\phi_i(\mathfrak{u}){\mathfrak{e}}^*_i(h)\Big)^{\!n}\!\!
=\exp\big(\phi_h(\mathfrak{u})\big).
\end{align*}
Similarly, applying $\varPsi$ to ${E}(h',h)$ in variable ${h'\in{H}}$, we  obtain
\begin{align*}
(\varPsi\circ{E})(\mathfrak{u},h)&=\Big|\sum_{n\ge0}\frac{1}{n!}\sum_{\imath^\lambda\vdash n}
\frac{n!}{\lambda!}\phi^\lambda_\imath(\mathfrak{u}){\mathfrak{e}}^{*\lambda}_\imath(h)\Big|^2
\Big(\sum_{n\ge0}\frac{1}{n!}\sum_{\imath^\lambda\vdash n}\frac{n!}{\lambda!}
|\mathfrak{e}^{*\lambda}_\imath(h)|^2\Big)^{-1}\\
&=\exp\big(2\mathop{\mathrm{Re}}\phi_h(\mathfrak{u})-\|h\|^2\big).
\end{align*}
By Lemma~\ref{Paley--Wiener},
$(\varPsi\circ{K})(\cdot,h)$ and $(\varPsi\circ{E})(\cdot,h)$ with ${h\in{H}}$ take values in
$ L_\chi^2 $. The expansion~\eqref{cp} follows from \cite[n.12]{Ito1952}
where polynomials $\mathfrak{h}_{n,m}(z,\bar{z})$ were introduced.
\end{proof}

\begin{theorem}\label{laplace1}
For any $f={\sum f_n\in L^2_\chi }$ with $f_n\in L^{2,n}_\chi$ the entire function
\[
\hat{f}(h):={\left\langle\exp(h)\mid\varPhi^* f\right\rangle_\beta}\quad\text{of variable}\quad{h\in H}
\]
and its Taylor coefficients at zero $d^n_0\hat{f}$ have the integral representations
\begin{equation}\label{laplaceA}\begin{split}
\hat{f}(h)&=\int \exp(\bar\phi_h)f\,d\chi=
\int \exp\big(2\mathop{\mathrm{Re}}\phi_h-\|h\|^2\big)f\,d\chi,\\
d^n_0\hat{f}(h)&=\int\bar\phi_h^nf_n\,d\chi,
\end{split}
\end{equation}
respectively. The Fourier transform ${F\colon L^2_\chi \ni{f}\longmapsto\hat{f}\in{H}^2_\beta}$ provides the isometries
 \[
 {L^2_\chi \stackrel{F}\simeq{H}^2_\beta}\quad\text{and}\quad
 L^{2,n}_\chi\stackrel{F}\simeq {P}_\beta^n(H).
 \]
\end{theorem}

\begin{proof}
Since $\varPsi=\varPhi\circ*^{-1}$, we obtain $\varPsi^*=*\circ \varPhi^*$.  From  \eqref{cformula} it follows that
$\hat{f}(h)=\left\langle\exp(h)\mid\varPhi^*f\right\rangle_\beta=
{\big\langle(\varPsi^*\circ f)(\cdot)\mid{K}(\cdot,h)\big\rangle_{{H}^2_\beta}}
=\big\langle(\varPsi^*\circ f)(\cdot)\mid{E}(\cdot,h)\big\rangle_{{H}^2_\beta}$. Thus,
\begin{align*}
\hat{f}(h)&={\big\langle(\varPsi^*\circ f)(\cdot)\mid{K}(\cdot,h)\big\rangle_{{H}^2_\beta}}=
{\big\langle(\varPsi^*\circ f)(\cdot)\mid{E}(\cdot,h)\big\rangle_{{H}^2_\beta}}\\
&={\big\langle f(\cdot)\mid(\varPsi\circ E)(\cdot,h)\big\rangle_\chi}
={\int \exp\big(2\mathop{\mathrm{Re}}\phi_h-\|h\|^2_H\big)f\,d\chi}
\end{align*}
by Lemma~\ref{infty2}. On the other hand, according to the same claim
\[
\hat{f}(h)=\big\langle(\varPsi^*\circ f)(\cdot)\mid{K}(\cdot,h)\big\rangle_{{H}^2_\beta}=
\big\langle f(\cdot)\mid(\varPsi\circ K)(\cdot,h)\big\rangle_\chi=\int \exp\left(\bar\phi_h\right)f\,d\chi.
\]

It particularly follows that for all $h=\alpha{x}$ with $x\in H$,
\[
\hat{f}\left(\alpha{x}\right)=
\int \exp\left(\bar\phi_{\alpha x}\right)f\,d\chi=
\sum\alpha^n\!\int \frac{\bar\phi_{x}^n}{n!}f_n\,d\chi, \quad{\alpha\in\mathbb{C}}.
\]
Using the $n$-homogeneity of derivatives, we find
\[
d^n_0\hat{f}(\alpha{x})=\frac{d^n}{d\alpha^n}\!
\sum\alpha^n\int \frac{\bar\phi_{x}^n}{n!}f_n\,d\chi
\mid_{\alpha=0}=\int\bar\phi_x^nf_n\,d\chi.
\]

Finally, we notice that the isometry ${L^2_\chi \stackrel{F}\simeq{H}^2_\beta}$ holds, since the isometry $\varPhi^*$ is surjective by Lemma~\ref{pbasis}.
Similarly, we get $L^{2,n}_\chi\stackrel{F}\simeq {P}_\beta^n(H)$.
\end{proof}

 \begin{corollary}\label{ggauss}
 For any $h\in{H}$ the Paley--Wiener map $\phi_h$ satisfies the equality
 \begin{equation*}
 \int\exp\big\{\mathop{\mathrm{Re}}\phi_h\big\}\,d\chi=
 \exp\Big\{\frac{1}{4}\|h\|^2\Big\}.
 \end{equation*}
\end{corollary}
\begin{proof}
It is enough to put $f\equiv1$ and to replace $h$ by $h/2$ in the formula \eqref{laplaceA}.
\end{proof}

 \begin{corollary}\label{PWF}
 The isometry ${*\colon\Gamma_\beta(H)\longrightarrow{H}^2_\beta}$ has the  factorization
$*={F}\circ\varPhi$.
 \end{corollary}
\begin{proof}
In fact,
${\varPhi\colon\Gamma_\beta(H)\ni\psi\longmapsto\varPhi\psi=f\in L^2_\chi }$ and
${F}\colon L^2_\chi \ni f\longmapsto
\hat{f}\in{H}^2_\beta$.
\end{proof}

\begin{corollary}\label{laplace1}
For every $f\in L^2_\chi $ the Taylor expansion at zero of the function
 \[
 \hat{f}(h)=\sum\frac{1}{n!}
 {d}^n_0\hat{f}(h)\quad\text{with}\quad
 f={\sum f_n\in L^{2}_\chi},\quad{f_n\in  L^{2,n}_\chi}
 \]
 has the coefficients
\begin{equation}\label{laplaceB}
d^n_0\hat{f}(h)=\int f_n\bar\phi_h^n\,d\chi=
\sum_{\imath^\lambda\vdash n}\hbar_\lambda
{s}^\lambda_\imath[{f}_\imath\,{\mathfrak{e}}_\imath^*(h)],\quad
{f}_\imath:=\int{f}\bar{\phi}_\imath\,d\chi
\end{equation}
with summation over all  standard Young tabloids ${[\imath^\lambda]}$
such that ${\imath^\lambda\vdash n}$
where $s^\lambda_\imath=0$ if the conjugate partition $\lambda^\intercal$ has $\lambda^\intercal_1>\eta(\lambda)$ and
${s}^\lambda_\imath[{f}_\imath\,{\mathfrak{e}}_\imath^*(h)]:=
s^\lambda_\imath(t_\imath)$
with $t_\imath={f}_\imath\,{\mathfrak{e}}_\imath^*(h)$.
\end{corollary}

\begin{proof}
By the Frobenius formula \cite[I.7]{MacDonald98} we find that
$\phi_h^n(\mathfrak{u})=
\sum_{\imath^\lambda\vdash n}\hbar_\lambda
s^\lambda_\imath(\mathfrak{u},h)$,
where $s^\lambda_\imath=0$ if  $\lambda^\intercal_1>\eta(\lambda)$, and
$s^\lambda_\imath(\mathfrak{u},h)$ is defined by \eqref{schur0}, whereas $\hbar_\lambda$ by \eqref{hook}.
Thus,
\begin{equation}\label{exp1}
\exp\phi_h(\mathfrak{u})=\sum_{n\ge0}\frac{1}{n!}
\sum_{\imath^\lambda\vdash n}\hbar_\lambda
s^\lambda_\imath(\mathfrak{u},h)=
\sum_{n\ge0}\frac{1}{n!}
\sum_{\imath^\lambda\vdash n}\frac{n!}{\lambda !}
\phi^\lambda_\imath(\mathfrak{u})\mathfrak{e}^{*\lambda}_\imath(h).
\end{equation}
Using  \eqref{exp1} in combination with  Theorem~\ref{Schur1}, we find
\begin{equation*}\label{FTayl}
\hat{f}(h)=\int f(\mathfrak{u})\exp\bar\phi_h(\mathfrak{u})\,d\chi(\mathfrak{u})
=\sum_{n\ge0}\frac{1}{n!}\sum_{\imath^\lambda\vdash n}
\hbar_\lambda\bar{s}^\lambda_\imath[{f}_\imath\,{\mathfrak{e}}_\imath^*(h)]
\end{equation*}
where the derivative at zero may be defined as
\[
d^n_0\hat{f}(h)=\sum\limits_{\imath^\lambda\vdash n}
\hbar_\lambda{s}^\lambda_\imath[{f}_\imath\,{\mathfrak{e}}_\imath^*(h)]\quad\text{with}\quad
{s}^\lambda_\imath[{f}_\imath\,{\mathfrak{e}}_\imath^*(h)]:=
\int f(\mathfrak{u})\bar{s}^\lambda_\imath(\mathfrak{u},h)\,d\chi(\mathfrak{u}).
\]
In fact, for $zh$ with  ${z\in\mathbb{C}}$ and $\imath^\lambda\vdash n$
with $\lambda^\intercal_1>\eta(\lambda)$
we find  \[{s}^\lambda_\imath[{f}_\imath\,\mathfrak{e}^*_\imath(zh)]=z^n
{s}^\lambda_\imath[{f}_\imath\,\mathfrak{e}^*_\imath(h)].\] Hence, the derivative
$d^n_0\hat{f}(h)=(d^n/dz^n)\hat{f}(zh)|_{z=0}$ is a Taylor coefficient of $\hat{f}$.

Now, the Frobenius formula and Theorem~\ref{Schur1} yield the first equality in \eqref{laplaceB}.
By Lemmas~\ref{pbasis} and \ref{infty2} the second formula in \eqref{laplaceB} also holds.
\end{proof}

\begin{remark}\label{cor7}
In the finite-dimensional case $\mathfrak{U}=U(m)$, the Hardy space ${H}^2_\beta$ of entire analytic
functions of variable $h\in\mathbb{C}^m$ has the following orthogonal basis
$\big\{\mathfrak{e}^{*\lambda}=\mathfrak{e}^{*\lambda_1}_1\ldots\mathfrak{e}^{*\lambda_m}_m
\colon{\lambda=(\lambda_1,\ldots,\lambda_m)\in\mathbb{Y}}\big\}$. The Fourier transform
\[
\hat{f}(h) =\int \exp(\bar\phi_h)f\,d\chi_m=
\int \exp\big(2\mathop{\mathrm{Re}}\phi_h-\|h\|^2\big)f\,d\chi_m,\quad h\in\mathbb{C}^m
\]
provides the surjective isometry
${F\colon L_{\chi_m}^2\ni{f}\longmapsto\hat{f}\in{H}^2_\beta}$, defined by mappings
\[
{F\colon\mathfrak{e}^{*\lambda}\mapsto\phi^\lambda}\quad\text{such that}\quad
\|\mathfrak{e}^{*\lambda}\|_{H^2_\beta}^2=
\|{\phi}^\lambda\|_{\chi_m}^2=\frac{(m-1)!\lambda!}{(m-1+|\lambda|)!}
\]
where the space $L_{\chi_m}^2$ with the Haar measure $\chi_m$ on $U(m)$ has  the orthogonal basis $\big\{\phi^\lambda=
{\phi}^{\lambda_1}_1\circ\pi_m^{-1}\ldots{\phi}^{\lambda_m}_m\circ\pi_m^{-1}\colon\lambda\in\mathbb{Y}\big\}$.
\end{remark}

\section{Intertwining properties of Fourier transform}\label{7}

The shift group on ${H}^2_\beta$ is defined as
$T_a\psi^*(h):=\left\langle\mathcal{T}_a\exp(h)\mid\psi\right\rangle_\beta$ for all $\psi\in\Gamma_\beta(H)$, ${a,h\in H}$.
By \eqref{cformula}, $\left\langle\mathcal{T}_a\exp(h)\mid\psi\right\rangle_\beta
=T_a\psi^*(h)=\big\langle T_a\psi^*(\cdot)\mid\exp\langle\cdot\mid h\rangle\big\rangle_{H^2_\beta}$. Hence,
\[
T_a\psi^*(h)=\left\langle\mathcal{T}_a\exp(h)\mid\psi\right\rangle_\beta\!=
\left\langle\psi^*(\cdot)\mid\exp\langle\cdot\mid h+a\rangle\right\rangle_{H^2_\beta}\!=
\left\langle\psi^*(\cdot)\mid{M}_{a^*}\exp\langle\cdot\mid h\rangle\right\rangle_{H^2_\beta}
\]
where ${M}_{a^*}\exp\langle\cdot\mid h\rangle:=\exp{a}^*(\cdot)\exp\langle\cdot
\mid h\rangle=\exp\langle\cdot\mid h+a\rangle$ is defined to be
the multiplicative group onto the total set
$\{{\exp\langle\cdot\mid h\rangle}\colon h\in H\}$ in ${H}^2_\beta$.

Comparing the above formulas, we obtain that ${M}_{a^*}$ is adjoint to ${T}_a$ on ${H}^2_\beta$.
By virtue of adjoint relations, $\|{T}_a\psi^*\|_{H^2_\beta}
=\|{M}_{a^*}\psi^*\|_{H^2_\beta}$. The isometry $H^2_\beta\simeq\Gamma_\beta(H)$
yields $\|{T}_a\psi^*\|_{H^2_\beta}=\|\mathcal{T}_a\psi\|_\beta$.
According to  \eqref{contract}, we have
\begin{equation}\label{contract1}
\begin{split}
& \|T_a\psi^*\|^2_{H^2_\beta}\le\exp\big(\|a\|^2\big)
 \|\psi^*\|^2_{H^2_\beta} \quad \text{and} \quad
{T}_{a+b}={T}_a{T}_b={T}_b{T}_a
\\
& \|M_{a^*}\psi^*\|^2_{H^2_\beta}\le\exp\big(\|a\|^2\big)
 \|\psi^*\|^2_{H^2_\beta} \quad \text{and} \quad
{M}_{a^*+b^*}={M}_{a^*}{M}_{b^*}={M}_{b^*}{M}_{a^*}
\end{split}
 \end{equation}
 for ${a,b\in H}$.
Thus, these groups are strongly continuous with densely defined closed generators
$\partial^*_a\psi^*:={\lim_{z\to0}(T_{za}\psi^*-\psi^*)/z}$ and
${a}^*\psi^*:={\lim_{z\to0}(M_{za^*}\psi^*-\psi^*)/z}$.

Hence, the additive group $(H,+)$ on ${H}^2_\beta$ is represented  by
${M}_{a^{*}}\colon{H}^2_\beta\rightarrow{H}^2_\beta$ and
the generator $dM_{za^{*}}/dz\mid_{z=0}=a^{*}$ of its $1$-parameter subgroup $M_{za^{*}}$ is
strongly continuous  with the dense domain $\mathfrak{D}(a^{\!*})={\big\{\psi^*\in{H}^2_\beta\colon
a^{*}\psi^*\in{H}^2_\beta\big\}}$. On the other hand, the group $(H,+)$ can be  represented  as
$M_{a^{*}}^\dagger=\varPsi{M}_{a^{*}}\varPsi^*
\colon L^2_\chi \rightarrow L^2_\chi $.
The generator of its strongly continuous subgroup
\[
{\mathbb{C}\ni z\longmapsto M^\dagger_{z{a^{*}}}},\quad dM^\dagger_{z{a^{*}}}/dz\mid_{z=0}=\bar\phi_a
\quad\text{with}\quad\bar\phi_a=\varPsi{a^{*}}\varPsi^*
\]
has the dense domain $\mathfrak{D}(\bar\phi_a)={\big\{f\in L^2_\chi \colon
\bar\phi_a f\in L^2_\chi \big\}}$ and is closed, since $a^*$ is closed.

The group $(H,+)$ on $L^2_\chi $ can be also represented  by
$T_a^\dagger:=\varPsi{T}_a\varPsi^*\colon L^2_\chi \rightarrow L^2_\chi$.
From Lemmas~\ref{T_a},  \ref{pbasis} it follows that the generator of
strongly continuous  subgroup
\[
{\mathbb{C}\ni z\longmapsto T^\dagger_{z\mathfrak{a}}},\quad dT^\dagger_{za}/dz\mid_{z=0}=\partial_a^\dagger
\quad\text{with}\quad\partial_a^\dagger:=\varPsi\partial_a^*\varPsi^*
\]
has the dense domain $\mathfrak{D}(\partial_a^\dagger)={\big\{f\in L^2_\chi \colon
\partial_a^\dagger f\in L^2_\chi \big\}}$ and is closed, since $\partial_a^*$ is closed.
By \eqref{cformula}
$\hat{f}(h)=\left\langle\exp(h)\mid\varPhi^*f\right\rangle_\beta=
\big\langle(\varPsi^*\circ f)(\cdot)\mid\exp\langle\cdot
\mid h\rangle\big\rangle_{H^2_\beta}$. Hence, by Lemma~\ref{infty2},
\begin{align*}
T_a^\dagger\hat{f}(h)=\big\langle(\varPsi^*\circ f)(\cdot)\mid T_a\exp\langle\cdot
\mid h\rangle\big\rangle_{H^2_\beta}=\int f\exp\left(\bar\phi_{h+a}\right)\,d\chi.
\end{align*}

\begin{lemma}\label{conn}\label{comm}
The additive group $(H,+)$ on $L^2_\chi $ has two representations
${a\mapsto M_{a^{*}}^\dagger}$ and ${a\mapsto T_a^\dagger}$
which are adjoint,  strongly continuous with closed densely defined
generators $\bar\phi_a$ and $\partial_a^\dagger$, respectively.
For every $f\in\mathfrak{D}(\bar\phi_a^m)={\big\{f\in L^2_\chi \colon
\bar\phi_a^m{f}\in L^2_\chi \big\}}$ with $m\in\mathbb{N}_0$,
\begin{equation}\label{conec}
\partial_a^{*m}T_a{F}(f)=
{F}\big(\bar\phi_a^mM^\dagger_{a^{*}}f\big),\quad a\in{H}.
\end{equation}
For every $f\in\mathfrak{D}(\partial_a^{\dagger m})={\big\{f\in L^2_\chi \colon
\partial_a^{\dagger m}{f}\in L^2_\chi \big\}}$ with $m\in\mathbb{N}_0$,
\begin{equation}\label{ai}
a^{* m}M_{a^{*}}{F}(f)={F}\big(\partial_a^{\dagger m}T_a^\dagger f\big),\quad{a}\in{H}.
\end{equation}
As a conclusion, $\partial_{\mathbbm{i}a}^\dagger=-\mathbbm{i}\partial_a^\dagger$. Moreover,
the following commutation relations hold,
\begin{equation}\label{comm2}
 M_{a^*}^\dagger T_b^\dagger=
\exp\langle{a}\mid{b}\rangle T_b^\dagger M_{a^*}^\dagger,
\qquad\big(\bar\phi_{a}\partial_b^\dagger-
\partial_b^\dagger\bar\phi_a\big)f={\langle{a}\mid{b}\rangle}f,
\end{equation}
for all  $f$ from the dense subspace $\mathfrak{D}(\bar\phi_a^2)\cap
\mathfrak{D}(\partial_b^{\dagger 2})\subset{L}^2_\chi$ and  nonzero  $a,b\in{H}$.
\end{lemma}

\begin{proof}
Using that ${T}_a$ and ${M}_{a^{\!*}}$ are adjoint, we find that
\begin{equation*}
\partial_a^{*m}T_a\hat{f}(h)
=\int \frac{d^mM_{za^{*}}^\dagger f}{dz^m}\Big|_{z=0}\exp\bar\phi_h\,d\chi
=\int (\bar\phi_a^mf)\exp\bar\phi_h\,d\chi,\quad m\ge0
\end{equation*}
for all ${f\in L^2_\chi }$.  This gives \eqref{conec}.
Since
$M_{a^{*}}\psi^*(h)=\big\langle\psi^*(\cdot)\mid{M}_{a^{*}}\exp\langle\cdot\mid h\rangle\big\rangle_{H^2_\beta}=\exp{a^{*}}(h)\,\psi^*(h)$, we obtain
\begin{equation}\label{M0}
\begin{split}
a^{* m}M_{a^{*}}\hat{f}(h)&=
\frac{d^mM_{za^{*}}\hat{f}(h)}{dz^m}\Big|_{z=0}
=\int \frac{d^mT_{za}^\dagger f}{dz^m}\Big|_{z=0}\exp\bar\phi_h\,d\chi\\
&=\int (\partial_a^{\dagger m}f)\exp\bar\phi_h\,d\chi\quad\text{with}\quad
f\in\mathfrak{D}(\partial_a^{\dagger m}),\quad \psi^*=\varPsi^*f.
\end{split}\end{equation}
This together with the group property by applying  ${F}$ and ${F}^{-1}$ yields \eqref{ai}.

Now, we prove the commutation relations. For any ${f\in L^2_\chi }$ and $h\in{H}$, we have
\begin{align*}
M_{b^{*}}{T}_a\hat{f}(h)&=
\exp\langle{h}\mid{b}\rangle\hat{f}(h+a),\\
T_a{M}_{b^{*}}\hat{f}(h)&=\exp\langle{h+a}\mid{b}\rangle\hat{f}(h+a)
=\exp\langle{a}\mid{b}\rangle{M}_{b^{*}}{T}_a\hat{f}(h).
\end{align*}

For each
$\hat{f}\in\mathfrak{D}(b^{*2})\cap\mathfrak{D}(\partial_a^{2})$
and ${t\in\mathbb{C}}$ by differentiation, we obtain
\begin{equation}\label{d2t2}
\big(d^2/dt^2\big)T_{ta}M_{tb^{*}}\hat{f}\mid_{t=0}=
\big(\partial_a^{*2}+2\partial_a^*b^{*}+b^{*2}\big)\hat{f}.
\end{equation}
Subsequently, taking into account \eqref{d2t2} together with
${(d/dt)[\exp\langle ta\mid \bar{t}b\rangle M_{tb^{*}}T_{ta}]}$
$=[(d/dt)\exp\langle{ta}\mid \bar{t}b\rangle]
M_{tb^{*}}T_{ta}+
\exp\langle{ta}\mid {\bar{t}b}\rangle
[(d/dt)M_{tb^{*}}T_{ta}],$ we find
\begin{align*}
\big(\partial_a^{*2}+2\partial_a^*b^{*}+b^{*2}\big)\hat{f}&
=(d/dt)\big[(d/dt)\exp\langle ta\mid\bar{t}b\rangle
M_{tb^{*}}T_{ta}\hat{f}\big]_{t=0}\\
&=2{ \langle{a}\mid{b}\rangle}\hat{f}+
\big(\partial_a^{*2}+2b^{*}\partial_a^*+b^{*2}\big)\hat{f}.
\end{align*}
Hence, for each  $\hat{f}$ from the dense subspace
$\mathfrak{D}(b^{*2})\cap\mathfrak{D}(\partial_a^{2})
\subset{H}^2_\beta$, which includes  all polynomials  generated by
finite sums $\varPsi^*(f)=\bigoplus\psi_n\in\Gamma_\beta(H)$ with
$\psi_n\in{H}^{\odot{n}}_\beta$,
\begin{align}\label{comm1}
T_aM_{b^{*}}&=
\exp\langle{a}\mid{b}\rangle M_{b^{*}}T_a,\quad
\left(\partial_a^*b^{*}-b^{*}\partial_a^*\right)\hat{f}={ \langle{a}\mid{b}\rangle}\hat{f}.
\end{align}

Corollary~\ref{PWF} yields
${F}=*\circ\varPhi^*$ and  ${F}^{-1}=\varPhi\circ*^{-1}$.
The  equality  \eqref{M0} for $m=0$ can be rewritten as
$M_{b^{*}}\hat{f}(a)=
\left\langle\exp(a)\mid {T}_b\varPhi^* f\right\rangle_\beta$
with ${f\in L^2_\chi }$ or in another way $*\circ{T}_b=M_{b^{*}}\circ*$.
Hence,
$T_b^\dagger=\varPhi\,{T}_b\varPhi^*
=\varPhi\circ{*}^{-1}\circ{M}_{b^*}\circ*\circ\varPhi^*=
{F}^{-1}M_{b^{*}}\,{F}$ and
$\partial_b^\dagger={F}^{-1}b^{*}\,{F}$. Similarly,
$M_{a^*}^\dagger ={F}^{-1}T_a\,{F}$
and $\bar\phi_a={F}^{-1}\partial_a^*\,{F}$. Finally,
\begin{align*}
&M_{a^*}^\dagger T_b^\dagger={F}^{-1}T_aM_{b^{*}}\,{F}
=\exp\langle{a}\mid{b}\rangle{F}^{-1}M_{b^{*}}T_a\,{F}
=\exp\langle{a}\mid{b}\rangle
 T_b^\dagger M_{a^*}^\dagger,\\
&\big(\bar\phi_a\partial_b^\dagger-\partial_b^\dagger\bar\phi_a\big)f=
{F}^{-1}\left(\partial_a^*b^{*}-b^{*}\partial_a^*\right){F}f=\langle{a}\mid{b}\rangle {f}
\end{align*}
for all $f$ from the dense subspace
$\mathfrak{D}(\bar\phi_a^2)\cap\mathfrak{D}(\partial_b^{\dagger 2})\subset{L}^2_\chi$, which includes  all functions  generated by
finite sums $\varPhi\left(\bigoplus\psi_n\right)$ with $\psi_n\in{H}^{\odot{n}}_\beta$.
\end{proof}

\section{Infinite-dimensional Heisenberg group}\label{8}
Our goal is to describe an irreducible representation on the space $L^2_\chi $
of the group $\mathcal{H}_\mathbb{C}$, defined by \eqref{HC}.
We will use the appropriate generalization of Weyl's system which in our case is written in the form of
$L^2_\chi $-valued function of variable ${h\in H}$
\[
{W}^\dagger(h):={W}^\dagger(a,b)=
\exp\Big\{\frac{1}{2}\langle{a}\mid{b}\rangle\Big\}T^\dagger_bM^\dagger_{a^*}.
\]

For convenience, we will use the quaternion algebra
$\mathbbm{H}=\mathbb{C}\oplus\mathbb{C}\mathbbm{j}$  of numbers
$\zeta={{(\alpha_1+\alpha_2\mathbbm{i})}+{(\alpha'_1+\alpha'_2\mathbbm{i})}\mathbbm{j}}
={\alpha+\alpha'\mathbbm{j}}$ such that
$\mathbbm{i}^2=\mathbbm{j}^2=\mathbbm{k}^2=
\mathbbm{i}\mathbbm{j}\mathbbm{k}=-1$,
$\mathbbm{k}=\mathbbm{i}\mathbbm{j}=-\mathbbm{j}\mathbbm{i}$,
$\mathbbm{k}\mathbbm{i} = -\mathbbm{i}\mathbbm{k}= \mathbbm{j}$,
where $(\alpha,\alpha')\in\mathbb{C}^2$ with
$\alpha={\alpha_1+\alpha_2\mathbbm{i}},\alpha'={\alpha'_1+\alpha'_2\mathbbm{i}\in\mathbb{C}}$ and
$\alpha_\imath,\alpha'_\imath\in\mathbb{R}$ $(\imath=1,2)$  \cite[5.5.2]{procesi2007lie}.
Let us denote $\alpha':=\Im{\zeta}$ for all $\zeta=\alpha+\alpha'\mathbbm{j}\in\mathbbm{H}$.

Consider  the Hilbert  space ${H}\oplus{H}\mathbbm{j}$
 with   $\mathbbm{H}$-valued inner product
\begin{align*}
\langle{h}\mid{h}'\rangle&=\langle{a}+{b}\mathbbm{j}\mid{a}'+{b}'\mathbbm{j}\rangle
=\langle{a}\mid{a}'\rangle+\langle{b}\mid{b}'\rangle+
\left[\langle{a}'\mid{b}\rangle-\langle{a}\mid{b}'\rangle\right]\mathbbm{j}
\end{align*}
where  ${h}={a}+{b}\mathbbm{j}$ with ${a}$, ${b}\in{H}$.
Hence,
\[
\Im\langle{h}\mid{h}'\rangle=\langle{a}'\mid{b}\rangle-\langle{a}\mid{b}'\rangle,\qquad
\Im\langle{h}\mid{h}\rangle=0.
\]

\begin{theorem}\label{heis}
The representation  of \,$\mathcal{H}_\mathbbm{C}$ over  $L^2_\chi $ in the Weyl-Schr\"{o}dinger form
\[
S^\dagger\colon\mathcal{H}_\mathbb{C}\ni X(a,b,t)\longmapsto
\exp(t){W}^\dagger(h),\quad{h}=a+b\mathbbm{j}
\]
 is well defined and irreducible. The  Weyl system satisfies the relation
\begin{equation}\label{weylS}
{W}^\dagger(h+h')=\exp\Big\{-\frac{\Im\langle{h}\mid{h}'\rangle}{2}\Big\}
{W}^\dagger(h){W}^\dagger(h')
\end{equation}
which on any real subspace $\{\tau{h}\colon \tau\in\mathbb{R}\}$ transforms to the $1$-parameter group
\begin{equation}\label{weylG}
{W}^\dagger\left((\tau+\tau'){h}\right)
={W}^\dagger(\tau{h}){W}^\dagger(\tau'{h})=
{W}^\dagger(\tau'{h}){W}(\tau{h})
\end{equation}
with the densely defined generator on  $L^2_\chi $  of the form $\mathfrak{p}^\dagger_{h}:=\partial_{b}^\dagger+\bar\phi_{a}$.
Moreover, the following commutation relations hold,
\begin{align}\label{Wcomm}
{W}^\dagger({h}){W}^\dagger({h}')&=
\exp\big\{\Im\left\langle{h}\mid{h}'\right\rangle \big\}
{W}^\dagger({h}'){W}^\dagger({h})\quad\text{where}\\\nonumber
\Im\left\langle{h}\mid{h}'\right\rangle&=-
\big[\mathfrak{p}_{h}^\dagger,\mathfrak{p}_{h'}^\dagger\big]
\quad\text{with}\quad
\big[\mathfrak{p}_{h}^\dagger,\mathfrak{p}_{h'}^\dagger\big]:=
\mathfrak{p}_{h}^\dagger\mathfrak{p}_{h'}^\dagger-
\mathfrak{p}_{{h}'}^\dagger\mathfrak{p}_{h}^\dagger
\end{align}
on the dense subspace $\mathfrak{D}(\bar\phi_a^2)\cap
\mathfrak{D}(\partial_b^{\dagger 2})\subset{L}^2_\chi$.
\end{theorem}

{\em Proof}. Let us
consider the auxiliary group $\mathbb{C}\times (H\oplus{H}\mathbbm{j})$ with multiplication
$(t,{h})(t',{h}')=
\left(t+t'-\frac{1}{2}\Im\langle{h}\mid{h}'\rangle,\,{h}+{h}'\right)$
for all ${h}={a}+{b}\mathbbm{j}$,
${h}'={a}'+{b}'\mathbbm{j}\in{H} \oplus{H} \mathbbm{j}$.
The mapping
${G}\colon
X({a},{b},t)\longmapsto
\left(t-\frac{1}{2}\langle{a}\mid{b}\rangle , \,{a}+{b}\mathbbm{j}\right)$
is a group isomorphism, since
\begin{align*}
&{G}\big(X({a},{b},t)X({a}',{b}',t')\big)
={G}\big(X({a}+{a}',{b}+{b}',
t+t'+\langle{a}\mid{b}'\rangle )\big)\\
&=\Big(t+t'+\langle{a}\mid{b}'\rangle
-\frac{1}{2}\big(
\langle{a}+{a}'\mid{b}+{b}'\rangle \big),
({a}+{a}')+
({b}+{b}')\mathbbm{j}\Big)\\
&=\Big(t+t'-\frac{1}{2}\big(\langle{a}\mid{b}\rangle +
\langle{a}'\mid{b}'\rangle \big)
+\frac{1}{2}\big(\langle{a}\mid{b}'\rangle -
\langle{a}'\mid{b}\rangle \big),({a}+{a})+
({b}+{b}')\mathbbm{j}\Big)\\
&=\Big(t-\frac{1}{2}\langle{a}\mid{b}\rangle , \,
{a}+{b}\mathbbm{j}\Big)
\Big(t'-\frac{1}{2}\langle{a}'\mid{b}'\rangle , \,
{a}'+{b}'\mathbbm{j}\Big)={G}\left(X({a},{b},t)\right)
{G}\left(X({a}',{b}',t')\right).
\end{align*}

On the other hand,  let us define the auxiliary Weyl system
\begin{equation}\label{Weyl0}
{W}({h})=
\exp\Big\{\frac{1}{2}\langle{a}\mid{b}\rangle \Big\}
M_{{b}^*}T_{{a}},
\quad{h}={a}+{b}\mathbbm{j}.
\end{equation}
Using group properties and the commutation relation \eqref{comm1}, we obtain
\begin{align}\nonumber
&\exp\Big\{-\frac{\Im\langle{h}\mid{h}'\rangle }{2}\Big\}
{W}({h}){W}({h}')
=\exp\Big\{\frac{\langle{a}\mid{b}'\rangle }{2}
-\frac{\langle{a}'\mid{b}\rangle }{2}\Big\}
{W}({h}){W}({h}')\\\label{Wcommut}
&=\exp\Big\{\frac{\langle{a}\mid{b}\rangle }{2}
+\frac{\langle{a}'\mid{b}'\rangle }{2}\Big\}
\exp\Big\{\frac{\langle{a}\mid{b}'\rangle }{2}
-\frac{\langle{a}'\mid{b}\rangle }{2}\Big\}
M_{{b}^*}T_{a}M_{{{b}'}^*}T_{{a}'}\\\nonumber
&=\exp\Big\{\frac{1}{2}\langle{a}
+{a}'\mid{b}+{b}'\rangle \Big\}
M_{{b}^*+{{b}'}^*}T_{{a}+{a}'}
={W}({h}+{h}').
\end{align}
Hence,  the  mapping ${\mathbb{C}\times ({H} \oplus{H} \mathbbm{j})\ni(t,{h})\longmapsto\exp(t){W}({h})}$
acts as a group isomorphism into the operator algebra  over ${H}^2_\beta$. So, the  representation
\[
S\colon\mathcal{H}_\mathbbm{C} \ni X({a},{b},t)\longmapsto
\exp(t){W}({h})=\exp\Big\{t+\frac{1}{2}\langle{a}\mid{b}\rangle \Big\}
M_{{b}^*}T_{{a}}
\]
 is also well defined over ${H}^2_\beta$, as a composition of  group isomorphisms.

Let us check the irreducibility. Suppose the contrary. Assume there exist an element ${h}_0\noteq 0$ in
${H} $ and an integer ${n>0}$ such that
\[
\exp\Big\{t+\frac{1}{2}\langle{a}\mid{b}\rangle \Big\}
\exp{\langle{c}\mid{a}\rangle }
\langle{c}+{b}\mid{h}_0\rangle ^n=0\quad\text{for all}\quad
{a},{b},{c}\in {H} .
\]
But, this is only possible for ${h}_0=0$. It gives a contradiction. Finally, using that
\[
\exp\Big\{t+\frac{1}{2}\langle{a}\mid{b}\rangle \Big\}
T_{b}^\dagger M_{a^*}^\dagger=
{F}^{-1}\Big(\exp\Big\{t+\frac{1}{2}\langle{a}\mid{b}\rangle \Big\}
M_{b^*}T_a\Big){F},
\]
we obtain that  $S^\dagger={F}^{-1}S\,{F}$ is   irreducible. Applying  ${F}$, ${F}^{-1}$ to \eqref{Wcommut} we get \eqref{weylS}.

Consider  the Weyl system ${W}^\dagger$ on the space $L^2_\chi $. By \eqref{weylS} we obtain the equality
\begin{align*}
&{W}^\dagger(h){W}^\dagger(h')=
\exp\Big\{\frac{\Im\langle{h}\mid{h}'\rangle }{2}\Big\}
{W}^\dagger({h}+{h}')=\exp\Big\{-\frac{\Im\langle{h}'\mid{h}\rangle }{2}\Big\}{W}^\dagger({h}'+{h})\\
&=\exp\Big\{\!\!-{\Im\langle{h}'\mid{h}\rangle }\Big\}
\exp\Big\{\frac{\Im\langle{h}'\mid{h}\rangle }{2}\Big\}{W}^\dagger({h}'+{h})
=\exp\Big\{\!\!-{\Im\langle{h}'\mid{h}\rangle }\Big\}
{W}^\dagger(h'){W}^\dagger({h}).
\end{align*}
Using  this equality,  we get \eqref{weylG} for any fixed ${h}={a}+{b}\mathbbm{j}\in{H}\oplus{H}\mathbbm{j}$.
The $1$-parameter group ${W}^\dagger(\tau{a},\tau{b})={W}^\dagger(\tau{h})$ with real $\tau$ has the  generator $\mathfrak{p}_{h}^\dagger=\mathfrak{p}_{a,b}^\dagger$, since
\[
\mathfrak{p}_{{a},{b}}^\dagger=\frac{d}{d\tau}{W}^\dagger(\tau{h})\Big|_{\tau=0}
=\frac{d}{d\tau}\exp\Big\{\frac{1}{2}\langle{\tau{a}}\mid {\tau{b}}\rangle \Big\}
T^\dagger_{\tau{b}}M^\dagger_{\tau{a}^*}\Big|_{\tau=0}=\partial_{b}^\dagger+\bar\phi_{a}.
\]
Taking into account the inequalities \eqref{contract1} and that $F$ is isometric, we get
 \[
 \|{W}^\dagger(\tau{a},\tau{b})f\|^2_\chi\le\exp\big(\|\tau a\|^2+\|\tau b\|^2\big)\|f\|_\chi^2,\quad f\in L^2_\chi.
 \]
 Hence, the group ${W}^\dagger(\tau{a},\tau{b})$ in variable $\tau\in\mathbb{R}$
 is strongly continuous on $L^2_\chi$ and therefore
has the dense domain $\mathfrak{D}(\mathfrak{p}_{h}^\dagger)
={\big\{f\in L^2_\chi \colon\mathfrak{p}^\dagger_{h}f\in L^2_\chi \big\}}$. Moreover, its generator  $\mathfrak{p}^\dagger_{h}$
is closed (see, e.g., \cite{Vrabie03}). Note also that $\mathfrak{p}^\dagger_{\tau{h}}=\tau\mathfrak{p}^\dagger_{h}$ for $\tau\in\mathbb{R}$.

Finally, applying the commutation relation \eqref{comm2} and commutability of group generators in different directions
over the dense set $\mathfrak{D}(\bar\phi_a^2)\cap\mathfrak{D}(\partial_b^{\dagger 2})\subset{L}^2_\chi$, we have
\begin{align*}
-\Im\langle{h}\mid{h}'\rangle &
=\left\langle{a}\mid{b}'\right\rangle -
\left\langle{a}'\mid{b}\right\rangle
=\bar\phi_{a}\partial_{{b}'}^\dagger-
\bar\phi_{{a}'}\partial_{b}^\dagger+
\partial_{b}^\dagger\bar\phi_{{a}'}-
\partial_{{b}'}^\dagger\bar\phi_{a}\\
&=(\partial_{b}^\dagger+\bar\phi_{a})
(\partial_{{b}'}^\dagger+\bar\phi_{{a}'})-
(\partial_{{b}'}^\dagger+\bar\phi_{{a}'})
(\partial_{b}^\dagger+\bar\phi_{a})=
\big[\mathfrak{p}_{h}^\dagger,\mathfrak{p}_{{h}'}^\dagger\big].
\end{align*}

\section{Heat equation associated with Weyl system}\label{9}

In what follows, we will consider the real Banach space $c_0$ and let $\xi_n^*$ be the coordinate functional, i.e., $\xi_n^*(\xi)=\xi_n$ for ${\xi\in c_0}$. Since, the embedding $\mathcal{I}\colon{l}_2\looparrowright c_0$ is continuous, the  Gelfand triple $l_1 \stackrel{\mathcal{I}^*}\longrightarrow l_2 \looparrowright c_0$ with adjoint $\mathcal{I}^*$ holds. The mapping $Q\colon{l}_1\to c_0$ with $Q:=\mathcal{I}\circ\mathcal{I}^*$
is positive and
$\langle Q\xi^*\mid Q\xi^*\rangle_{l_2}:=\xi^*(Q\xi^*)=\sum \xi_n^2=\|\xi\|^2_{l_2}$
where  $\xi=Q\xi^*\in\mathscr{R}(Q)$ and $\xi^*\in{l}_1=c^*_0$.
By the Aronszajn–Kolmogorov decomposition theorem (see e.g., \cite[Prop.1]{Weron1981}) the
appropriative reproducing kernel Hilbert space can be determined as  $\overline{\mathscr{R}(Q)}=l_2$.

Consider the abstract Wiener space defined by $\mathcal{I}\colon l_2\looparrowright c_0$.
Given $\xi_1^*,\ldots,\xi_n^*\in l^1=c_0^*$, we assign the family of cylinder sets
$\Omega_n^c=\left\{\xi\in c_0\colon(\xi_1^*(\xi),\ldots,\xi_n^*(\xi))\in\Omega_n\right\}$
with any Borel $\Omega_n\subset\mathbb{R}^n$ that are not a $\sigma$-field.
Define the $\sigma$-additive extension $\mathfrak{w}$ of the Gaussian  measure $\gamma$
onto the Borel $\sigma$-algebra $\mathscr{B}(c_0)$, called futhure the Wiener measure, such that
\[\mathfrak{w}(\Omega_n^c):=\gamma(\Omega_n)\quad\text{with}\quad
\gamma(\Omega_n):=(2\pi)^{-n/2}\int_{\Omega_n}
\exp\big\{-\|\omega\|_{l_2}^2/2\big\}\,d\omega.\]

By Gross' theorem \cite{Gross65} there exists a smaller abstract Wiener space $\{w_0,\|\cdot\|_{w_0}\}$
such that injections ${l_2\looparrowright{w}_0\looparrowright c_0}$ are continuous
and the increasing sequence of orthogonal projectors $p_n\colon{l}_2\to\mathbb{R}^n$  has the extension
$(p^\sim_n)$ on $w_0$ that is convergent to the identity operator on $w_0$ and ${\mathfrak{w}(w_0)=1}$.
The integral of any cylinder function ${\upsilon\colon c_0\to\mathbb{R}}$ such that
$\upsilon=\rho\circ p_n^\sim$ is defined to be
$\int_{\Omega_n^c}\upsilon\,d\mathfrak{w}=\int_{\Omega_n}\rho\,d\gamma$.
The Fernique theorem \cite{Fernique70},\cite[Thm~3.1]{Kuo75} implies that  these exist
$\varepsilon,\eta>0$ such that $\|\cdot\|_{w_0}$ satisfies the following conditions
with a sufficiently large $K>0$,
\[
\int_{c_0}\exp\big\{\varepsilon\|\xi\|^2_{w_0}\big\}d\mathfrak{w}(\xi)<\infty,\quad
\mathfrak{w}\big(\|\xi\|_{w_0}\ge K\big)\le\exp\big\{-\eta K^2\big\}.
\]

Let us go back to the Weyl system ${W}^\dagger$. Consider in $L^2_\chi$ the dense subspace
${L}^{+2}_\chi:=\bigcup_{n\ge0}\bigoplus_{m=0}^nL^{2,m}_\chi$.
Let  ${a}={b}=\mathbbm{i}\xi_m\mathfrak{e}_m$  with ${\xi_m\in\mathbb{R}}$. Then
by Theorem~\ref{heis}
\[
{W}^\dagger(\mathbbm{i}\xi_m\mathfrak{e}_m,\mathbbm{i}\xi_m\mathfrak{e}_m)
= \exp\big\{{-\xi_m^2}/{2}\big\}
T^\dagger_{\mathbbm{i}\xi\mathfrak{e}_m}M^\dagger_{-\mathbbm{i}\xi\mathfrak{e}_m^*}.
\]

\begin{theorem}\label{schrod}
For any $f\in{L}^{+2}_\chi$ and $\xi=(\xi_m)\in c_0$ there exists the limit
\[
{W}^\dagger_\xi{f}=\lim_{n\to\infty}{W}^\dagger_{p_n^\sim(\xi)}{f},\quad
{W}^\dagger_{p_n^\sim(\xi)}:=
\exp\Big\{-\frac{\|{p_n^\sim(\xi)}\|_{w_0}^2}{2}\Big\}
\prod_{m=1}^{n} T^\dagger_{\mathbbm{i}\xi_m\mathfrak{e}_m}
M^\dagger_{-\mathbbm{i}\xi_m\mathfrak{e}_m^*}
\]
$\mathfrak{w}$-almost everywhere on $c_0$ such that the $1$-parameter Gaussian semigroup
\begin{equation}\label{Gauss}
\mathfrak{G}^\dagger_rf=\frac{1}{\sqrt{4\pi r}}\int_{c_0}
\exp\Big\{-\frac{\|\xi\|_{w_0}^2}{4r}\Big\}{W}^\dagger_\xi{f}\,d\mathfrak{w}(\xi),
\quad r>0
\end{equation}
on the space ${L}^{+2}_\chi$ is  generated by $-\sum\big(\partial_m^\dagger+\bar\phi_m\big)^2$
with $\partial_m^\dagger:=\partial_{\mathfrak{e}_m}^\dagger$.
As a consequence, $w(r)=\mathfrak{G}^\dagger_rf$ is unique solution of the Cauchy problem
\begin{equation}\label{Hamilt}
\frac{dw(r)}{dr}=-\sum\big(\partial_m^\dagger+\bar\phi_m\big)^2w(r),
\quad w(0)=f\in{L}^{+2}_\chi.
\end{equation}
\end{theorem}

\begin{proof}
Note that $(M_{b^*}T_a)^*=T_a^*M_{b^*}^*=M_{a^*}T_b$. Hence,
$(\partial_{a}^\dagger+\bar\phi_{a})^*=\partial_{a}^\dagger+\bar\phi_{a}$ is self-adjoint  for $a=b$, as a
generator of the group
${W}^\dagger(\tau{a},\tau{a})=\exp\left\{{\|\tau a\|^2}/{2}\right\}T_{\tau{a}}^\dagger M_{\tau{a}^*}^\dagger$ with ${\tau\in \mathbb{R}}$.
Replacing $a=b$ by $\mathbbm{i}\tau a$ with ${\tau\in\mathbb{R}}$, we obtain that
\[
{W}^\dagger(\mathbbm{i}\tau{a},\mathbbm{i}\tau{a})=
\exp\Big\{-\frac{1}{2}\langle{\tau{a}}\mid {\tau{a}}\rangle \Big\}
T^\dagger_{\mathbbm{i}\tau{a}}M^\dagger_{-\mathbbm{i}\tau{a}^*} \quad
\text{has the generator} \quad \mathbbm{i}(\partial_{a}^\dagger+\bar\phi_{a})
\]
with self-adjoint $\partial_{a}^\dagger+\bar\phi_{a}$. By relations \eqref{comm2},
${W}^\dagger(\mathbbm{i}\tau{a},\mathbbm{i}\tau{a})$ is unitary.

Lemma~\ref{comm} implies that $[M_{-\mathbbm{i}\xi_m\mathfrak{e}_m^*}^\dagger,
T_{\mathbbm{i}\xi_k\mathfrak{e}_k}^\dagger]=0$ and
${[M_{-\mathbbm{i}\xi_m\mathfrak{e}_m^*}^\dagger, M_{-\mathbbm{i}\xi_k\mathfrak{e}_k^*}^\dagger]=0}$,
as well as,
$[T_{\mathbbm{i}\xi_m\mathfrak{e}_m}^\dagger, T_{\mathbbm{i}\xi_k\mathfrak{e}_k}^\dagger]=0$
for any $m\noteq k$. In view of the relations \eqref{comm2},
\begin{equation}\label{symp}
\big[\bar\phi_{\mathbbm{i}\xi_m\mathfrak{e}_m},
\partial_{\mathbbm{i}\xi_k\mathfrak{e}_k}^\dagger\big]=0 \quad \text{if}\quad m\noteq k
\quad\text{and}\quad
\big[\bar\phi_{\mathbbm{i}\xi_m\mathfrak{e}_m},
\partial_{\mathbbm{i}\xi_m\mathfrak{e}_m}^\dagger\big]=-\xi^2_m.
\end{equation}

Check that \eqref{Gauss} holds. Denote
${W}^\dagger_{p_n^\sim(\xi)}:=
\prod_{m=1}^{n}{W}^\dagger(\mathbbm{i}\xi_m\mathfrak{e}_m,\mathbbm{i}\xi_m\mathfrak{e}_m)$ and
$T^\dagger_{p_n^\sim(\xi)}:=\prod_{m=1}^{n} T^\dagger_{\mathbbm{i}\xi_m\mathfrak{e}_m}$, as well as,
$M^\dagger_{p_n^\sim(\xi)}:=\prod_{m=1}^{n} M^\dagger_{-\mathbbm{i}\xi_m\mathfrak{e}_m^*}$ with
${\xi=(\xi_m)\in{w}_0}$.  Using  \eqref{contract1} with the operator norm over $H^2_\beta$, we get
the inequality
\[
\ln\prod_{m=1}^{n}\|T_{\mathbbm{i}\xi_m\mathfrak{e}_m}\|_{\mathscr{L}(H^2_\beta)}^2
\le\sum_{m=1}^{n}\langle\xi_m\mathfrak{e}_m\mid\xi_m\mathfrak{e}_m\rangle^2
=\sum_{m=1}^{n}\xi_m^2=\|p_n^\sim(\xi)\|_{l_2}^2.
\]
The relation $T_{\mathbbm{i}\xi_m\mathfrak{e}_m}^\dagger
=\varPsi{T}_{\mathbbm{i}\xi_m\mathfrak{e}_m}\varPsi^*$ implies that the left-hand side term above can be changed by
$\ln\prod_{m=1}^{n}\|T_{\mathbbm{i}\xi_m\mathfrak{e}_m}^\dagger\|_{\mathscr{L}({L}^2_\chi)}^2$.
For $M^\dagger_{p_n^\sim(\xi)}=\prod_{m=1}^{n} M^\dagger_{-\mathbbm{i}\xi_m\mathfrak{e}_m^*}$ similarly.

Using the unitarity of groups ${W}^\dagger(\mathbbm{i}\xi_m\mathfrak{e}_m,\mathbbm{i}\xi_m\mathfrak{e}_m)$,
we find by virtue of \eqref{symp} that their product  ${W}^\dagger_{p_n^\sim(\xi)}=
\exp\left\{-\|{p_n^\sim(\xi)}\|_{l_2}^2/2\right\}T^\dagger_{p_n^\sim(\xi)} M^\dagger_{p_n^\sim(\xi)}$ is also unitary.
Taking into account the continuity of ${\mathcal{I}_0\colon l_2\looparrowright{w}_0}$ and
that $p_n^\sim$ converges  to the identity mapping on  $w_0$, as well as, that ${\mathfrak{w}(w_0)=1}$,
we obtain for all ${f\in{L}^{+2}_\chi}$, $n\ge0$,
\[
\|{W}^\dagger_{p_n^\sim(\xi)}f\|_\chi\le
\exp\big\{-\|{p_n^\sim(\xi)}\|_{l_2}^2/2\big\}\|{f}\|_\chi\le
\exp\big\{-\|\mathcal{I}_0\|^2\,\|\xi\|_{w_0}^2/2\big\}\|{f}\|_\chi.
\]

The Lebesgue dominated convergence theorem implies that there exists
$\lim\|{W}^\dagger_{p_n^\sim(\xi)}f\|_\chi$ $\mathfrak{w}$-almost everywhere in variable ${\xi\in{w}_0}$ for all ${f\in L^{2,m}_\chi}$ and ${m>0}$.
By completeness of $L^{2,m}_\chi$, the limit ${W}^\dagger_\xi{f}$ is well defined $\mathfrak{w}$-almost everywhere and
\begin{equation}\label{est}
\|{W}^\dagger_\xi{f}\|_\chi\le\exp\big\{-\|\mathcal{I}_0\|^2\,
\|\xi\|_{w_0}^2/2\big\}\|{f}\|_\chi\quad\text{for all}
\quad{f\in{L}^{+2}_\chi},\quad{\xi\in{w}_0}.
\end{equation}

The $\|\cdot\|_\chi$-norm of integrant in \eqref{Gauss} is bounded
by $\exp\left\{\varepsilon\|\xi\|_{w_0}^2\right\}$ with any $\varepsilon>0$.
By Fernique’s theorem and \eqref{est}, the integral \eqref{Gauss} with the Wiener measure $\mathfrak{w}$
exists for all $f\in{L}^{+2}_\chi$. The equality $\mathfrak{w}(w_0)=1$ implies that
the integral \eqref{Gauss} is absolutely convergent uniformly in variables $r>0$ on the whole space $c_0$.
It provides the $C_0$-property of $\mathfrak{G}_r$ in variables $r>0$ on any finite sum
$\bigoplus_{m=0}^nL^{2,m}_\chi$.

Prove that the semigroup
$\mathfrak{G}_r$ is generated by  $\sum\mathfrak{p}_m^{\dagger2}$
with $\mathfrak{p}^\dagger_m:={\mathbbm{i}(\partial_m^\dagger+\bar\phi_m)}$.
By differentiation of
${W}^\dagger(\mathbbm{i}\xi_m{a},\mathbbm{i}\xi_m{a})$ at ${\xi_m=0}$,
we get that its generator coincides with $\mathfrak{p}^\dagger_m$. In fact,
${W}^\dagger(\mathbbm{i}\xi_m{a},\mathbbm{i}\xi_m{a}){f}=\exp\big\{\xi_m\mathfrak{p}_m^\dagger\big\}f$ for all ${f\in\phi^\mathbb{Y}}$.
Applying the next  formula for Gamma functions with $\alpha=(\alpha_1,\ldots,\alpha_n)\in\mathbb{N}_0^n$
\[\begin{split}
\prod_{m=1}^{n}\!\frac{1}{\sqrt{4\pi r}}
\int\exp\left\{\frac{-\xi_m^2}{4r}\right\}\xi_m^{2\alpha_m}d\xi_m
\Big|_{\xi_m=2\sqrt{r}x_m}\!\!=\!\prod_{m=1}^{n}\!\frac{(2\sqrt{r})^{2}}{\sqrt{\pi}}\!
\int\exp\big\{\!-x_m^2\big\}x_m^{{2\alpha_m}}dx_m\\
={2^{2n}r^n}\!\prod_{m=1}^{n}\!\Gamma\Big(\frac{{2\alpha_m}+1}{2}\Big)
=2^nr^n\frac{({2\alpha}-1)!}{(\alpha-1)!},
\end{split}\]
we find that for any ${L}^{+2}_\chi$-valued cylinder function
$h_n=({W}^\dagger_\xi{f})\circ p_n^\sim$ we have
\begin{align*}
\mathfrak{G}^\dagger_rh_n&=\prod_{m=1}^{n}\frac{1}{\sqrt{4\pi r}}\int
\exp\Big\{-\frac{\xi^2_m}{4r}\Big\}
\exp\big\{\xi_m\mathfrak{p}_m^\dagger\big\}d\xi_mh_n\\
&=\sum_{\alpha\in\mathbb{N}_0^n}
\prod_{m=1}^{n}\frac{\mathfrak{p}_m^{\dagger\alpha_m}}{\alpha_m!}
\frac{1}{\sqrt{4\pi r}}
\int\exp\Big\{-\frac{\xi_m^2}{4r}\Big\}\xi_m^{\alpha_m}\,d\xi_mh_n\\
&=\sum_{\alpha\in\mathbb{N}_0^n}2^nr^n\prod_{m=1}^{n}\dfrac{(2\alpha_m-1)!}{(\alpha_m-1)!}
\dfrac{\mathfrak{p}_m^{\dagger2}}{(2\alpha_m)!}h_n
=\exp\Big\{r\sum_{m=1}^{n}\mathfrak{p}_m^{\dagger 2}\Big\}h_n.
\end{align*}
Using \eqref{est}, we obtain that $0\le r\longmapsto \mathfrak{G}_r^\dagger$ is the $1$-parameter $C_0$-semigroup on any finite sum
$\bigoplus_{m=0}^nL^{2,m}_\chi$ with densely defined closed generator $\sum_{m=1}^{n}\mathfrak{p}_m^{\dagger 2}$.
Applying the known  relation  \cite{Vrabie03} between the initial  problem \eqref{Hamilt} and the
$1$-parameter $C_0$-semigroup $\mathfrak{G}_r^\dagger$,  we obtain that the function $w_n(r)=\mathfrak{G}^\dagger_rf_n$
for any ${n\in\mathbb{N}}$ solves  this problem in the sense that $d\mathfrak{G}^\dagger_rf_n/dr|_{r=0}=\sum_{m=1}^{n}\mathfrak{p}_m^{\dagger2}f_n$
for all ${f_n\in}\bigoplus_{m=0}^nL^{2,m}_\chi$. The theorem is proved.
\end{proof}

Taking into account the isometries
${H}^2_\beta\stackrel{\varPsi}\simeq{L}^2_\chi $ and
${P}_\beta^n(H)\stackrel{\varPsi}\simeq {L}^{2,n}_\chi $ from \eqref{PSI},
defined by linearization, we can rewrite the Cauchy problem in polynomial form.

Consider the Weyl system ${W(a,b)=\exp\left\{\langle{a}\mid{b}\rangle/2 \right\}
M_{b^*}T_a}$ defined by \eqref{Weyl0} on the dense subspace of polynomials
${P}_\beta(H):={\sum}_{n\ge0}{P}_\beta^n(H)$ in ${H}^2_\beta,$
consisting of all finite sums of $n$-homogenous polynomials $\psi^*(h)=\sum\psi_n^*(h)$ of variable ${h\in H}$ with components $\psi_n^*=\mathcal{P}\circ\psi_n\in {P}_\beta^n(H)$.
Replacing $a$ by $\tau a$ and $b$ by $\tau b$ with real ${\tau\in\mathbb{R}}$, we get that
$T_{\tau a}$ and $M_{\tau b^*}$ are  generated by closed generators on ${P}_\beta(H)$,
\[
\partial^*_a\psi^*=\lim_{\tau\to0}\left(T_{\tau a}\psi^*-\psi^*\right)/\tau\quad\text{and}\quad
{a}^*\psi^*=\lim_{\tau\to0}\left(M_{\tau a^*}\psi^*-\psi^*\right)/\tau, \quad{a,b\in{H}}.
\]
As a consequence, the $1$-parameter Weyl system ${W}(\tau{a},\tau{b})$ has the  generator
\[
\frac{d}{d\tau}{W}(\tau{a},\tau{b})|_{\tau=0}
=\frac{d}{d\tau}\exp\Big\{\frac{1}{2}\langle{a}\mid{b}\rangle \Big\}\Big|_{\tau=0}=
b^*+\partial_a^*
\]
densely defined on ${P}_\beta(H)$ such that
$(\tau b)^*+\partial_{\tau a}^*=\tau(b^*+\partial_a^*)$ for real $\tau$.
Let ${W}_{p^\sim_n(\xi)}=
\prod_{m=1}^{n}{W}(\mathbbm{i}\xi_m\mathfrak{e}_m,\mathbbm{i}\xi_m\mathfrak{e}_m)$,
$T_{p^\sim_n(\xi)}=\prod_{m=1}^{n} T_{\mathbbm{i}\xi_m\mathfrak{e}_m}$,
$M_{p^\sim_n(\xi)}=\prod_{m=1}^{n} M_{-\mathbbm{i}\xi_m\mathfrak{e}_m^*}$.

\begin{corollary}\label{schrod2}
For all $\psi^*\in{P}_\beta(H)$ and $\xi=(\xi_m)\in c_0$ there exists the limit
\[
{W}_\xi\psi^*=\lim_{n\to\infty}{W}_{p^\sim_n(\xi)}\psi^*,\quad
 {W}_{p^\sim_n(\xi)}:=
\exp\Big\{-\frac{\|{p^\sim_n(\xi)}\|_{w_0}^2}{2}\Big\}
\prod_{m=1}^{n} M_{-\mathbbm{i}\xi_m\mathfrak{e}_m^*}
T_{\mathbbm{i}\xi_m\mathfrak{e}_m}
\]
$\mathfrak{w}$-almost everywhere on $c_0$ such that the $1$-parameter Gaussian semigroup
\begin{equation*}\label{Gauss1}
\mathfrak{G}_r\psi^*=\frac{1}{\sqrt{4\pi r}}\int_{c_0}
\exp\Big\{\frac{-\|\xi\|_{w_0}^2}{4r}\Big\}{W}_\xi\psi^*d\mathfrak{w}(\xi),
\quad r>0
\end{equation*}
is  generated
by $-\sum(\mathfrak{e}_m^*+\partial_m^*)^2$. Thus, $w(r)=\mathfrak{G}_r\psi^*$ is unique solution
of the problem
\begin{equation*}\label{Hamilt1}
\frac{dw(r)}{dr}=-\sum\big(\mathfrak{e}_m^*+\partial_m^*\big)^2w(r),
\quad w(0)=\psi^*\in{P}_\beta(H)
\end{equation*}
in the space of Hilbert--Schmidt  polynomials ${P}_\beta(H)$.
\end{corollary}

\footnotesize

\bibliographystyle{spmpsci}      

\end{document}